\def\**{$*\!*\!*$}


 at 12truept
\font\tenmsy=msym10
\textfont8=\tenmsy
\mathchardef\ssm="7872

\input psfig
\mathsurround = 2pt
\abovedisplayskip=6pt
\belowdisplayskip=6pt

\def \QED {\rlap{$\sqcup$}$\sqcap$\smallskip}

\def\ref{\hangindent=1pc \hangafter=1 \noindent}

\input psfig


\centerline{\bf Combinatorics, geometry and attractors of 
 quasi-quadratic maps.}\medskip

\centerline{Mikhail Lyubich\footnote{*}{\rm Supported in part by 
NSF grant DMS-8920768 and a Sloan Research Fellowship.}}\smallskip
\centerline{Mathematics Department and IMS, SUNY Stony Brook}\bigskip
\centerline{February 1992, revised October 1992 }\bigskip

\noindent{\bf Abstract.} The Milnor problem on one-dimensional attractors
  is solved for $S$-unimodal maps with a non-degenerate critical point $c$.
  It provides us with a complete understanding of the possible
  limit behavior for Lebesgue almost every point. This theorem follows
   from a geometric study of the critical set $\omega(c)$ of a 
  ``non-renormalizable" map. It is proven
  that the scaling factors characterizing the geometry of this set
  go down to 0 at least exponentially. This resolves the problem of the
  non-linearity control in small scales. The proofs strongly involve
   ideas from  renormalization theory and holomorphic dynamics.
\bigskip

\centerline{\bf \S1. Introduction.}\smallskip

Let $f:[0,1]\rightarrow [0,1]$ be an $S$-unimodal map 
(see the definitions later)
of the interval with a non-degenerate critical point $c$.
Let us call such a map  {\sl quasi-quadratic}.
As usual, $\omega(x)$ denotes the limit set of the forward orb$(x)$.
The following theorem 
solves the Milnor problem [M1].

\proclaim Theorem on the Measure-Theoretic Attractor.
 Let $f$ be a quasi-quadratic map normalized
by the condition $f: c\mapsto 1\mapsto 0$.
Then there is a unique set $A$ (a measure-theoretic attractor in the 
sense of Milnor) such
that $A=\omega(x)$ for Lebesgue almost all $x\in [0,1]$, and only one
of the following three possibilities can occur:
{\item 1.} $A$ is a limit cycle;
{\item 2.} $A$ is a cycle of intervals;
{\item 3.} $A$ is a Feigenbaum-like attractor. 

This result gives a clear picture of measurable dynamics for the maps under
consideration. Let us explain the words used in the statement.
A limit cycle is the periodic orbit whose basin of attraction has 
non-empty interior. A cycle of intervals is the union of finitely many 
intervals $I_n$ with disjoint interiors cyclically interchanged 
by the dynamics. A Feigenbaum-like
attractor is an invariant Cantor set of the following structure:
$$ A=\bigcap{\cal O}_n,$$
where ${\cal O}_1\supset {\cal O}_2\supset...$ is a nested sequence of
cycles of intervals of increasing periods.

It makes sense to compare the above Theorem with its topological counterpart
known since late 70s:

\proclaim  Theorem on the Topological Attractor ([MT], [G], [JR], [vS]). 
Let $f$ be an $S$-unimodal map normalized
by the condition $f: c\mapsto 1\mapsto 0$.
Then there is a unique set $\Lambda$ (a topological attractor) such
that $\Lambda=\omega(x)$ for a generic  $x\in [0,1]$, and only one
of the following three possibilities can occur:
{\item (i).} $\Lambda$ is a limit cycle;
{\item (ii).} $\Lambda$ is a cycle of intervals;
{\item (iii).} $\Lambda$ is a Feigenbaum-like attractor.

From this point of view the Theorem on the Measure-Theoretic Attractor 
says that the map $f$ has a unique
measure-theoretic attractor $A$ coinciding with the topological
attractor $\Lambda$. In cases (i) and (iii) this
 was proven by Guckenheimer [G].
In case (ii) it is known that
$$\Lambda=\bigcup_{k=0}^{p-1} I_k$$
where $I_i$ form  a cycle of intervals of period $p$, and
$f^{\circ p}|\; I_k$ is {\sl topologically exact}. The last property means
that for any interval $J\subset I_k$ there is an $n$ such that
$f^{\circ pn} J=I_k$.
So, we can reduce the Theorem on the Measure-Theoretic Attractor
 to the following statement: 


\proclaim Theorem A.  Let $f$ be a quasi-quadratic map. Assume 
that $f$ is topologically exact on $[0,1]$, hence $\Lambda=[0,1]$. 
Then $\omega(x)=[0,1]$ for Lebesgue almost all $x\in [0,1]$,
that is $A=[0,1]$.

This measure-theoretic result  follows from geometric properties
of {\sl the critical set} $\omega(c)$.
 In order to study microstructure
of this set  we introduce the concept of a generalized
renormalization as the rescaled first return map restricted to a
neighborhood of the critical set. This moves us out of the class 
of unimodal maps to a class ${\cal T}$ of maps with a single critical point 
but defined on the union of disjoint intervals mapped onto a bigger
interval.\footnote{$^{\dag}$}{Up to some point we allow the domain to
 consist of
infinitely many intervals. However, the renormalization philosophy becomes
really valuable only in the finite case.}
 Any map with recurrent critical point becomes 
infinitely renormalizable in this sense.  

Let us consider the sequence of renormalized maps
$$R^{\circ n} f: \bigcup_i I^n_i\rightarrow I^{n-1}_0.$$
As basic geometric characteristics of the critical set  we consider the
 scaling factors $\mu_n=|I^n_0|/|I^{n-1}_0|$ 
and the Poincar\'{e} lengths of the gaps between the intervals $I^n_i$
of level $n$. Our main geometric result is

\proclaim Theorem B. Let $f$ be a topologically exact quasi-quadratic map
with the recurrent critical point.
Then the  Poincar\'{e} lengths of the gaps of 
the renormalized maps go up to $\infty$.

 Moreover, to a first approximation they increase at least linearly.
Unfortunately, there is one unpleasant
circumstance which can slow the rate down, 
namely the long cascades of ``central returns"
(when the map is combinatorially  close to being renormalizable).
However, the rate is still under explicit control. In particular,
if we number these cascades by $\kappa\in {\bf N}$ then the rate
 will be at least linear in $\kappa$.

As it is clear from the very name of our field,
the basic problem of non-linear dynamics is to gain control
of non-linearity. In our situation the non-linearity in small scales
is controlled by the above mentioned   scaling factors. Theorem B
 implies that
{\sl the scaling factors go down at least exponentially in} $\kappa$.
 This provides us with a
perfect control of non-linearity: the high order renormalized maps
are becoming purely quadratic exponentially fast. 
Observe  that this  contrasts drastically with the 
Feigenbaum-like case   when the geometry of 
the critical set is bounded from below 
(provided the   combinatorics is bounded)
 which creates a definite amount of
 non-linearity in all scales  (see Sullivan [S]). 

Let us now dwell in more detail on the ideas of this work.
 \S 2 contains the combinatorial treatment of maps of class ${\cal T}$. 
 According to the recurrent properties of the orb($c$) we 
split the analysis into two subcases: the {\sl reluctantly} and 
the {\sl persistently} recurrent. The latter case presents a stronger
recurrence of the critical point. For example, in this case the return
time of points of the orb$(c)$ back to a neighborhood $U\ni c$ is
bounded.
Surprisingly this helps to provide the
further analysis since the renormalized maps turn out to be of finite
type (that is, defined on the finitely many intervals).
 We show that such maps are 
classified by the cascades of renormalization types, and that these types
can be combined independently. The last statement gives us  a big freedom
in producing examples. In particular, the so called Fibonacci map naturally
arises as a map of the simplest stationary type. This map studied in
[LM] was  a basic model which  clarified the situation.



\S\S 3-6 are occupied with the proof of Theorem B. 
In   \S 3 we prove it under the assumption that we
start with a small enough scaling factor. 
To this end we introduce
 a modification of the Poincar\'{e} lengths,
the ``asymmetric Poincar\'{e} lengths",  which behave more regularly
under renormalization. We show that the asymmetric Poincar\'{e} lengths 
of the gaps increase almost monotonically. 
Though this is the longest technical piece of the paper, we are actually in
a quite comfortable position since the starting condition provides us at once
with a perfect control of non-linearity.  

In order to get rid of the starting assumption we pass to the complex plane 
in the class of {\sl (generalized) polynomial-like maps}, a complex
counterpart of maps of class ${\cal T}$ (\S 6). These maps
are  defined on the finite union of disjoint topological disks mapped
onto a bigger disk.
The crucial fact is that all such maps with the same combinatorics are 
quasi-conformally equivalent which yields quasi-symmetric equivalence
on the real line. This part relies on recent developments
in holomorphic dynamics, particularly on the work of Branner and Hubbard [BH].
Since the conclusion of Theorem B is quasi-symmetrically invariant, 
it is enough to produce just one example of a polynomial-like map $f$ with a
given combinatorics and arbitrarily small starting scaling factor.
Exactly at this point it is important to have a freedom  gained by
passing to the class of {\sl generalized} polynomial-like maps
and its real counterpart, class ${\cal T}$.

A bridge between quadratic maps and generalized polynomial-like maps
is given by the renormalization on a Yoccoz puzzle-piece introduced
in [L2], a complex analogue of the renormalization mentioned above. 
This passage completes the proof of Theorem B for quadratic polynomials
(see the end of \S 6).
The argument for general quasi-quadratic maps  
is based upon two {\sl Sullivan's Principles} [S]: 
{\item 1)} {\sl Renormalizations of sufficiently smooth maps are becoming
real analytic;}
{\item 2)} {\sl High order renormalization of a real analytic map makes 
it polynomial-like.}

\noindent
These principles turn out to be efficient in our setting of generalized 
renormalization.
The First Principle works when  we have combinatorics of bounded type, 
that is the number of intervals $I^n_i$ is bounded on all levels.
To push it forward we need a priori distortion bounds obtained by
Martens [Ma]. As to the case of unbounded combinatorics (as well as the
non-minimal case), it is actually
easier, and can be treated by a purely real argument. We discuss these
issues in \S 4.

In \S 5 we push the Second Principle forward. In order to construct a
polynomial-like map we consider the complex pull-backs of Euclidian disks based
upon the real intervals, and estimate their sizes (``complex bounds").  
This is another technical piece of the work. 

Finally, in the last \S 7 we derive Theorem A from Theorem B
by showing that the set $X=\{x: \omega(x)\ni c\}$ has
zero lower density at $c$.
This is our first application of geometric Theorem B but we expect 
a lot of others.

\medskip\noindent
{\bf  Remarks. 1.} Since  Milnor's paper of 1985 there has been a good deal
of effort to attack the problem of attractors, see
 [BL1-3], [GJ], [K], [Ma], [JS]. In particular, it was already known
 that  there is only one measure-theoretic attractor (see [BL1-3] or
 [GJ], [K]). 
 The reluctantly recurrent case was resolved in [BL3] and [GJ]. 
 
The geometric study of ``non-renormalizable maps" was started in [GJ],
 [Ma] and [JS]. In [JS] the absence of Cantor attractors was proved for
topologically exact  maps sufficiently close to the Chebyshev polynomial
$x\mapsto 4x(1-x)$.
 
 In [HK] the Fibonacci map was suggested as a candidate for a ``wild"
  situation when the  measure-theoretic and topological attractors are
  different. The paper [LM] showed that it is not the case for
 quasi-quadratic maps (an alternative purely real argument
  has been recently found in [KN]). However, a computer experiment 
  carried out jointly
  with F.Tangerman indicates that this may well be true in higher degrees.

  For a survey on the problem of attractors see [L3]. 

\smallskip {\bf 2.} The complex counterpart of Theorem A says that the
Lebesgue measure of the Julia set of a ``non-tunable"
 quadratic polynomial is equal to
0 (Lyubich [L2] and Shishikura (unpublished)). The case of cubic polynomials
with Cantor Julia set had been earlier treated by McMullen (see [BH]).
 The  real result turns out to be harder than the  complex one
because conformal invariants (like annulus moduli or Poincar\'{e} metric)
have only quasi-invariant analogues in the real setting. However, Theorem B
yields that on deep levels they are  becoming invariant exponentially fast
with respect to pull-backs. This makes the real-complex dictionary much
more precise.

\smallskip {\bf 3.} There are two pieces of the paper which strongly depend
on the quadratic-like nature of the critical point.
First, the growth of the Poincar\'{e} lengths of the gaps relies on the
fact that the square root map divides the Poincar\'{e} length at most
by 2 (Lemma 3.8). (By the way, this is the place where the asymmetric
 Poincar\'{e} length comes to the scene). The second place is the 
Branner-Hubbard divergence property (\S 6).

\medskip
\noindent {\bf Some definitions, notations and conventions.} A continuous map
$f: [0,1]\rightarrow [0,1]$ is called {\sl unimodal} if it has only
one extremum $c$. It is called $S$-{\sl unimodal} if additionally it is
three times differentiable with  only critical point $c$ and with negative
Schwarzian 

 \eject 

\noindent derivative outside $c$:
$$Sf={f'''\over f'}-{3\over 2}\left({f''\over f'}\right)^2<0.$$ 
If additionally the critical point $c$ is non-degenerate, that is 
$f''(c)\neq 0$, then the map is called {\sl quasi-quadratic}.
 
Let $\phi(x)=(x-c)^2$.
We use the abbreviations qs and qc for ``quasi-symmetric" and
``quasi-conformal" respectively. 

Saying ``an interval" we mean a
``closed interval". The notation $I=[a,b]$ means that $a$ and $b$ are the
endpoints of $I$ but does not necessarily mean that $a<b$. Points $a$ and $a'$
are called  ``$c$-symmetric" if $f a=f a'$. So, we can also talk about
$c$-symmetric sets.

The $n$-fold iterate of a map $g$ is  denoted by $g^{\circ n}$. The 
orb($x)\equiv g$-${\rm orb}(x)$ denotes $\{g^{\circ m} x\}_{m=0}^\infty$. Also
orb$_n(x)=\{g^{\circ m} x\}_{m=0}^n$ denotes the initial
piece of the orbit. Let ${\bf N}=\{0,1,2,...\}$.

We recommend the forthcoming book of de Melo and van Strien [MS]
 for the background in  one-dimensional dynamics.

\medskip\noindent
{\bf Acknowledgement.}
I'd like to thank Marco Martens for useful discussions, and for reading 
parts of the manuscript.
I owe to John Milnor several nice suggestions which improved the 
exposition.  I also take this opportunity to thank IHES and IMPA
for their hospitality while I was doing parts of this work.

\bigskip

\centerline{\bf \S2. Renormalization.}\medskip

Renormalization in dynamical systems means the first return map to
an appropriate piece of the phase space and then rescaling of this
piece to the ``original size". For example, if we have a periodic
$c$-symmetric interval $I\subset [0,1]$ of period $p>1$, we can consider
the unimodal map $f^p|I$ and then rescale $I$  back to [0,1]. This is the
usual notion of  renormalization in one-dimensional dynamics.
So, we can talk about renormalizable and non-renormalizable maps.
Repeating this procedure we can further talk about ``twice renormalizable",
``thrice renormalizable",..., or at the end ``infinitely renormalizable"
 maps. For example, looking through the Theorem on the Topological Attractor
we see that case (ii) corresponds to at most finitely renormalizable maps,
while the last case (iii) corresponds to infinitely renormalizable maps.
The case we concentrated on, when the map $f$ is topologically exact 
on $[fc, f^2 c]$, is equivalent to being non-renormalizable.

However, this specific terminology highly accepted in one-dimensional dynamics
is somewhat misleading. Indeed, we will see that the most interesting
``non-renormalizable" maps can  actually be treated as {\sl infinitely
renormalizable} in an appropriate sense, and this gives an  efficient
tool for studying the geometry of the critical set. This renormalization
does not respect the class of unimodal maps. So, let us describe an
appropriate class of maps.  

\medskip\noindent{\bf Class}${\cal T}.$   
Let $q,p\in [0,\infty]$, and let $I_k,\; k=-q,...,p,$  be a finite family
of disjoint intervals  compactly contained in a {\sl base} interval $J$.
The interval $I_0$ is marked, and will be called {\sl central}.
Let us consider a map 
$$g:\; {\rm Dom}(g)=\bigcup_{k=-q}^p I_k\rightarrow J \eqno (2-1)$$
 with a single turning point $c\in I_0$ (we will still call it ``critical"),
 which also
satisfies the following property: $g (\partial I_k)\subset \partial J$. 
In particular,
$g$ maps  each non-central interval $I_k,\; k\neq 0$, homeomorphically onto
the whole interval $J$, while the central interval $I_0$ is $c$-symmetric.

Let us also assume that $g$ does not have non-trivial wandering intervals.
This condition holds automatically under some regularity assumptions,
e.g., $g$ has negative Schwarzian derivative, and the critical point $c$ is
non-flat (see [G] or [L1]).  

Denote this class of maps by ${\cal T}^{\infty}$. Let us say that 
$g\in {\cal T}^{\infty}$ is
{\sl of finite type} if its domain  consists of only finitely many
intervals, that is $p,q\in {\bf N}$. Let ${\cal T}$ denote the subclass
of maps of finite type. Unimodal maps can also be viewed
as maps of class ${\cal T}$ whose domain contain a single interval.

A  point $x$ is said to be {\sl non-escaping} if
 $g^{\circ n} x\in {\rm Dom}(g)$ for all $n=1,2,...$.
 Denote by $K(g)$ the set of
non-escaping points ( the ``filled-in Julia set" of $g$).



\medskip\noindent {\bf Pull-Backs and nice intervals}.
The following {\sl pull-back construction} plays an essential role in 
what follows. Given an interval $T$ and  a point $x$ such that 
$g^{\circ n} x\in {\rm int}T $, we can pull $T$ back along the orbit 
$\{ g^{\circ k} x\}_{k=0}^n$. This means that we inductively construct a
sequence of intervals
$$T_n\equiv T\ni g^{\circ n} x,\quad T_{n-1}\ni g^{\circ(n-1)} x,
   \quad ...,\, T_0\ni x $$
so that $T_k$ is the maximal interval containing $g^k x$ whose image
is contained in $T_{k+1}$.  An interval $T_k$ of the pull-back is called
{\sl critical} if int$T_k\ni c$.  Since $g(\partial I_k)\subset \partial J$,

$$g(\partial T_k)\subset \partial T_{k+1},\; k=0,1,...,n-1.  \eqno (2-2)$$ 
Hence the non-critical intervals $T_k$ diffeomorphically map onto
$T_{k+1}$, and the critical intervals $T_k$ are $c$-symmetric,
$k=0,1,...,n-1$.
A pull-back is called {\sl monotone} if all intervals $T_k$
except  perhaps $T$ are non-critical. It is called {\sl unimodal}
if $T_0$ is critical while $T_1,...,T_{n-1}$ are non-critical. 

As in [Ma] let us call a $c$-symmetric interval $T$ {\sl nice} if
$$g^{\circ n}(\partial T)\cap T = \emptyset,\; n=1,2,...   \eqno (2-3)$$
For example, let $\alpha$ be a periodic point, and let $\beta$ be its
$n$-fold preimage, $\beta'$ be a $c$-symmetric point. Then the interval
$T=[\beta, \beta']$ is nice provided $\beta\not\in {\rm orb }(\alpha)$.
It follows that there are nice intervals in any neighborhood of $c$
(provided $g$ has no limit cycles).

Let us denote by ${\cal M}\equiv {\cal M}_T$ the family of
 intervals obtained by pulling
 a nice interval $T$ back along all possible orbits 
(${\cal M}$ comes from ``Markov"). Following the 
analogy with the holomorphic
setting (compare [H], [M2] or \S 6),
 the intervals of ${\cal M}$ will also be called {\sl puzzle-pieces}.
We use the notation $T^{(n)}(x)$ for the puzzle-piece obtained by pulling 
$T$ along orb$_n(x)$, and call $n$ the {\sl depth} of the puzzle-piece.
 (If we  start with another interval, say, $J$  then, of course,
we use the notations $J^{(n)}(x)$ for the corresponding puzzle-pieces).
The basic properties of the family  ${\cal M}$ are:\smallskip

\noindent
 (i) any two intervals of ${\cal M}$ are either disjoint or
  ``strongly nested";  in the latter case the interval of higher depth
    is contained in the interior  of the other one;\smallskip

\noindent
 (ii) any non-critical puzzle-piece $T^{(n)}(x)$ diffeomorphically
  maps onto $T^{(n-1)}(f x)$.\smallskip

\noindent
 (iii)  For any critical puzzle-piece  $T^{(n)}(x)$ of depth $>0$
     $$f (T^{(n)}(x), \partial T^{(n)}(x)) 
    \subset (T^{(n-1)}(f x), T^{(n-1)}(f x)).$$
    All critical puzzle-pieces are $c$-symmetric. \smallskip

\noindent (iv) If $c$ is recurrent then for any puzzle-piece
 $T^{(n)}(x)$, the orb$(c)$ does not cross $\partial T^{(n)}(x)$.\smallskip

\proclaim Lemma 2.1. Let $K\in{\cal M}$ be a critical puzzle-piece,
  and $x$ be a point whose
  orbit crosses {\rm int}$K$. Take the first moment $n>0$ for which
 $g^{\circ n} x\in {\rm int } K$. 
 Then the pull-back of $K$ along orb($x$) is either
 monotone or unimodal. 

\noindent{\bf Proof.} Let $K=K_m, K_{m-1},..., K_0\ni x$ be the above 
pull-back. Then the puzzle-pieces $K_{m-1}$,..., $K_0$ have  higher depths
 than $K$. It follows that they are are disjoint from $K$
and hence non-critical.             \QED

In particular, if $x\not\in K$ then the pull-back is monotone.
If $x=c$ is the critical point then the pull-back is certainly unimodal.

Let $K, L\in {\cal M}$ be two critical intervals, $L\subset K$. The interval
 $L$ is called a {\sl kid} of $K$ if it is obtained by a {\sl unimodal}
pull-back along a piece of the critical orbit 
$\{g^{\circ k} c\}_{k=0}^n,\; g^n c\in K$. The kid corresponding to the
first return of the critical point back to  $K$ will be called the 
{\sl first kid} $K^1$. Lemma 2.1 can be improved in the following way.

\proclaim Lemma 2.2. Let $K\in{\cal M}$ be a critical puzzle-piece, $L$
  be its first kid, and $x$ be  a point whose
  orbit crosses {\rm int}$L$. Take the first moment $n>0$ for which
 $g^{\circ n} x\in {\rm int } L$. Then the pull-back of $K$ along orb($x$) is 
 monotone if $x\not\in L$, and unimodal otherwise.

\noindent{\bf Proof.} Consider the subsequent
return moments $0<n_1<n_2<...<n_l=n$ of orb$(x)$ to $K$ until the first 
meeting with $L$.  According to the previous lemma, 
the first return produces a unimodal or monotone pull-back
of $K$ depending on whether $x\in L$ or not. All further returns produce 
monotone pull-backs of $K$. 
This is what was required.          \QED

\medskip\noindent {\bf Sequence of first grandkids. Cascades of central 
returns.}
Let $I^0\equiv J$ be the base interval, $I^1\equiv I_0$ be its first kid,
 $I^2$ be the first kid of $I^1$ etc.
In such a way we construct a nested sequence of the 
``first grandkids"
$$I^0\supset I^1\supset I^2\supset...$$
which will play fundamental role in what follows.
Let $t(n)$ be the first return time of the critical point back to $I^{n-1}$.
We say that the return to level $n-1$ is  {\sl central} if 
$g^{\circ t(n)}c\in I^n$.
Otherwise the return is classified as {\sl high} or {\sl low  }  depending on
whether $g^{\circ t(n)} I^{n}\supset I^n$ or 
$g^{\circ t(n)} I^{n}\cap I^n=\emptyset$
(compare [GJ]).

Denote by ${\cal L}\subset {\bf N}$ the sequence of all levels
  $m$ such that $g^{\circ t(m)} c\in I^{m-1}\ssm I^m$
 which means that the return to
the level $m-1$ is non-central. Let 
$\kappa: {\bf N}\rightarrow {\bf N}$ be the monotone surjective map
such that $\kappa(m+1)=\kappa(m)+1$ for $m\in {\cal L}$ and
$\kappa(m+1)=\kappa(m)$ otherwise.  The series of levels with the same
$\kappa$ will be called {\sl the cascades of central returns}.
(The ``cascade" can degenerate to a single level
if the central return does not actually occur).
 So, $\kappa$ numbers subsequently the cascades of central returns. 

An important issue  which we will discuss next is whether the sequence of
grandkids shrink down to the critical point.
( The main concern of this work will be
the {\sl rate of shrinking}).

\medskip\noindent 
{\bf Unimodal renormalization.}
  Let us say that $g\in {\cal T}^{\infty}$
 {\sl admits a unimodal renormalization}
if there exists a $c$-symmetric periodic interval. Then the return map to
this interval is unimodal which justifies the terminology. It will be convenient
to consider unimodal maps as ``admitting unimodal renormalizations".
So, when we say that a map $g\in {\cal T}^{\infty}$ does not admit unimodal
renormalizations, we assume automatically that it is not unimodal itself.
(The standard meaning of a renormalizable unimodal
map corresponds to admitting a unimodal renormalization with a period $>1$.)


\proclaim Lemma 2.3. The following properties are equivalent:
{\item (i)} $g$ does not admit a unimodal renormalization.
{\item (ii)} The sequence of levels with non-central returns is infinite.
{\item (iii)} The grandkids $I^n$  shrink to the critical point.
{\item (iv)} The  set $K(g)$ of non-escaping points has empty interior.

\noindent
{\bf Proof.} Let $t(n)$ be as above the return time of $c$ back to $I^{n-1}$.
Clearly, $t(n)$ is monotonically increasing.

\smallskip
(i) $\Rightarrow $ (ii). Assume that the returns to all
levels $m\geq n-1$ are central.  Then   $\cap_{m\geq n} I^m$ is a $c$-symmetric
$g^{\circ l}$-invariant interval. Contradiction. \smallskip

\smallskip
(ii) $ \Rightarrow $(iii).  If $t(m)=t,\; m\geq n,$ eventually stabilizers
then the returns to all levels $\geq n-1$ should be central contradicting (ii).
If $t(n)\to\infty$ but $I^n$ don't shrink to the critical point then
the intersection $\cap I^n$ is a wandering interval.

\smallskip (iii) $ \Rightarrow$ (iv). Since an appropriate iterate of $I^n$ covers
 the whole base interval $I^0=J$, there are  escaping points in $I^n$.
 Since $I^n$ shrink down to the critical point, there are escaping points
 in any neighborhood of $c$. But an appropriate iterate of any  other 
 interval must cover the critical point (no homtervals). Hence, escaping points
  are dense.

\smallskip (iv) $\Rightarrow$ (i).
  Any periodic interval is contained in the filled-in Julia set. 
\QED

The condition (iv) can also be stated in the following way:
$${\rm diam} (J^{(n)}(x))\to 0\quad {\rm as}\quad n\to\infty \eqno (2-4)$$
uniformly in $x$. If $g$ is of finite type then it follows that $K(g)$
is a Cantor set.
 
\smallskip\noindent {\bf Remarks. 1.}
One can also  easily see that a map $g\in {\cal T}$
admits a unimodal renormalization if and only if there is a an interval
$I^{n-1}$ such that the critical point does not escape its kid $I^{n}$ under
iterates of $f^{\circ t(n)}$. In the complex setting it will be accepted
as the main definition. \smallskip

\noindent {\bf 2.} One can also state and prove {\sl the Two kids Lemma}
by the same argument as in the complex setting  (see \S 6). However,
we don't need it for the real discussion.

\medskip\noindent {\bf Persistent and reluctant recurrence.}

Let $B(x,\epsilon)$ denote the interval centered at $x$ of radius $\epsilon$.

Assume  $c$ is recurrent. 
It is called {\sl reluctantly  recurrent} if
there exist an $\epsilon>0$ and an arbitrary long backward orbit 
$\bar x=\{x, x_{-1},... x_{-l}\}$ in $\omega(c)$  such that the $B(x,\epsilon)$
 allows a monotone  pull-back along $\bar x$. Otherwise $c$ is called
{\sl  persistently recurrent}.

\proclaim Lemma 2.4. Assume $g$ does not admit unimodal renormalizations.
 Then the following two properties are equivalent:
{\item (R1)} $c$ is reluctantly recurrent;\smallskip 
{\item (R2)} There is a critical puzzle-piece $I^{(m)}(c)$ 
              with infinitely many kids.

\noindent {\bf Proof.} $(R1)\Rightarrow (R2)$. By (2-4) we can find a 
puzzle-piece $J^{(n)}(x)\subset B(x,\epsilon)$. It is certainly a 
monotone pull-back of some critical puzzle-piece $J^{(m)}(c)$.

Let us pull  $J^{(n)}(x)$ to
depth $n+l$ along the backward orbit $\bar x$. Then let us consider the 
first moment when orb$(c)$ crosses $J^{(n+l)}(x_{-l})$ and pull 
this puzzle-piece further to the critical point. We get a kid of
$J^{(m)}(c)$. Letting $l\to \infty$ we obtain infinitely many kids. 

\smallskip $(R2)\Rightarrow (R1)$. Consider the pull-backs of
 $J^{(m)}(c)$  along the backward orbits of $c$ creating the kids. \QED 

\noindent {\bf Remark.} Persistently recurrent situation
 appeared in [BL3, Lemma 11.1]  and [GJ] under different names.
 In the complex setting it was introduced in [Y] under
 the name we use here. The term ``reluctantly recurrent" was suggested
  by McMullen. 

An invariant set $K$ is called {\sl minimal}
if the orbits of all points $x\in K$ are dense in $K$.

\proclaim Lemma 2.5 (see [BL3], [Ma]).
 In the persistently recurrent case the critical set $\omega(c)$ 
 is a minimal Cantor set.

\noindent{\bf Proof.} Assume there exists an $x\in\omega(c)$ whose orbit
does not accumulate on $c$. Then for sufficiently small $\epsilon>0$
the pull-backs of $B(g^{\circ n} x,\epsilon)$ along the orb$_n(x)$
are monotone. The contradiction proves that $\omega(c)$ is minimal. 
In particular, $\omega(c)$ does not contain periodic points. Hence it
is a Cantor set.   \QED

The minimality property yields
that for any relative neighborhood $U\subset\omega(c)$
all points $x\in U$ eventually return back to $U$, and, moreover, the
return time is bounded (but certainly depends on $U$). So, on $\omega(c)$
the return map  is {\sl of finite type}: it is defined on the 
 finitely many intervals.
This motivates the following consideration.

\medskip\noindent {\bf Return maps of finite type.}
Take a nice interval $L\ni c$. Let $K_0\ni c$ be its first kid.
Select also finitely many other pairwise disjoint pull-backs 
$K_i\subset L,\; -n\leq i\leq m$ corresponding to the first returns 
of some points back to $L$. 
Then we can define the first return map of class ${\cal T}$:
 $$h: \bigcup_{j=-n}^m K_j\rightarrow L. \eqno (2-5)$$
In the most interesting case when $L=I_0$ the combinatorial type of $h$
can be described in terms of the $g$-itineraries of $K_i$ through
the intervals $I_k$ of the previous level. To this end let us prepare
a bit of algebraic language.

\medskip\noindent {\bf Spin semigroups.}
Let $\Gamma$ be a semi-group. We call $\Gamma$
a {\sl spin  semi-group} if it is supplied with a character
$\epsilon: \Gamma\rightarrow {\bf Z}^2$. 
If $\Gamma$ is free then the spin structure can be certainly prescribed
arbitrary on the generators, and uniquely determined by this data.

Let $p,q\in {\bf N}$, and let us consider a free semi-group
$$\Gamma=<I_{-q},..., I_0,..., I_p> $$
with $p+q+1$ generators which are linearly ordered with a marked
generator $I_0$. 
If $\epsilon(I_0)>0$ we say that $\Gamma$ is {\sl positively oriented.}
 Otherwise $\Gamma$ is {\sl negatively oriented}. 

Let $W(\Gamma)\subset \Gamma$ be the set of words
 $\gamma=I_{k(0)}...I_{k(s)}$ such that   the last symbol is $I_0$, 
while all others are different from $I_0$. 
 Let us order $W(\Gamma)$ as follows (compare [MT]).
 Let $\gamma_t=I_{k(0)}...I_{k(t)}$ be the initial part of 
$\gamma$, $t=0,...,s$. Assign to $\gamma$ an integer
vector with components $\epsilon(\gamma_{t-1}) k(t)$, $t=1,...,s$  .
The ordering of $\Gamma$ is induced by  the lexicographic ordering of
the corresponding vectors. Note that any two words of $W(\Gamma)$
are comparable since they end with the $I_0$-symbol.

So, the object we have described should be called {\sl a free ordered spin
semigroup.} However, we will usually call it just ``spin semigroup"
keeping in mind the other structures. 

Let $(\Gamma', \epsilon')$ be another spin semigroup of the same type with
generators \break $I'_{-q'},...,I'_0,...,I'_{p'},$ and let
$\chi: (\Gamma',\epsilon')\rightarrow (\Gamma, \epsilon)$ be a semigroup
homomorphism. Let us call it  {\sl unimodal} if
$\chi(I_k)\in W(\Gamma)$ and 
$\chi$ is unimodal on the set of generators, that is,
 it is strictly monotone on the sets   $I'_{-q'},...,I'_0$ and
$I'_0,...,I'_{p'}$, and has an extremum at $I_0$. Moreover,
we require $I_0$  to be the minimum or maximum depending on whether $\Gamma$
is positively or negatively oriented.

A  unimodal homomorphism $\chi$ is called {\sl admissible} if additionally
 $$\epsilon'(K_j)={\rm sgn}j\; \epsilon(\chi(K_j).  \eqno (2-6)$$
\noindent

\medskip\noindent
{\bf  Return type.}
Take a map $g\in {\cal T}$ as in (2-1) of finite type.
 Let the numbering of the intervals $I_k$ be consistent
with the line order.
 Consider a free group $\Gamma_g$ generated by $I_k$.
Let us supply $\Gamma_g$ with a spin structure $\epsilon_g$. 
Spin of  a non-central interval $I_k$ 
 is equal to +1 or -1 
depending on whether $g|I_k$ preserves or reverses orientation. 
The spin of the central interval $I_0$  is equal to +1 or -1 depending on 
whether $c$ is the minimum or maximum point. 

Assume now that the critical point returns to $I_0$.  Then we can consider 
a return map (2-5):
$$h: \bigcup_{j=-n}^m K_j\rightarrow I_0$$
of class ${\cal T}$.  
Let $(G_h, \; \epsilon_h)$ be an associated spin semigroup, and 
let us define a homomorphism 
$$\chi:(\Gamma_h,\epsilon_h)\rightarrow (\Gamma_g,\epsilon_g)$$
assigning to each interval $K=K_j$  
its itinerary through the intervals of the previous level
$$\chi(K)=
 I_{k(1)}...I_{k(s-1)}\; I_0\in W(\Gamma_g) \eqno (2-7)$$ 
 until the first return back to $I_0$, that is, 
$g^j K\subset I_{k(j)}$ with $k(j)\neq 0$ for $j=1,...,s-1.$
 Let us call the itinerary $\chi(K_0)\in W(\Gamma_g)$ of the central
interval the {\sl kneading sequence} of the return map.


We call the homomorphism $\chi$ the {\sl type  of the return map}. 
Observe that $\chi$ is admissible. Property (2-6) expresses the chain rule
for the orientation type of the composition. The unimodal
property of $\chi$ reflects the fact that the central branch of $g$ is unimodal.



Denote by ${\cal T} (q,p,\epsilon)$ the subclass of ${\cal T}$ 
with specified $q,p\i {\bf N}$ and spin character $\epsilon$. 
Let $g^t\in{\cal T}(q,p,\epsilon),\; t_0\leq t\leq t_1$,
 be a continuous (in $C^0$ topology)  one-parameter family of maps.
This means that the intervals  $I_k\equiv I_k(t)$ 
move continuously with $t$, and each
rescaled branch $g_t|I_k(t)\circ\psi_k^t$ depends continuously
on $t$ (here $\psi_k^t: [0,1]\rightarrow I_k(t)$ is the 
orientation preserving rescaling).  Let us call the family $g^t$  {\sl full} if
$g^{t_0} I_0(t_0)$ covers all intervals $I_k(t_0)$ while $g^{t_1} I_0(t_1)$
 does not intersect the interiors of $I_k(t_1)$ (or vice versa).

\proclaim Lemma 2.6.
Let $g^t\in {\cal T}(q,p,\epsilon),\; 0\leq t \leq 1$,
 be a full one-parameter family.
Let $\chi: \Gamma'\rightarrow \Gamma_g$ be an admissible homomorphism.
Then for some  parameter value $t\in (0,1)$  there is a return map
$h^t\in {\cal T}$ whose return  type is equal to $\chi$. Moreover,
there is a full family of return maps
$h^t, \; t_0\leq t\leq t_1,$ with the same property.

\noindent {\bf Proof.} Without loss of generality we can assume that
the critical point is minimum for all $g^t$, that is, $\Gamma_g$ is
positively oriented. In what follows we omit $t$
in notations keeping in mind that everything actually depends on it.
For every $t$ let us consider
the intervals $L_j$ with itineraries $\chi(K_j)$. Such intervals
exist (because outside the central interval $g$ acts as the Bernoulli scheme),
and continuously depend on $t$.
Moreover, $g^{\circ(l_j-1)} L_j = I_0$ where $l_j$ is the length of the
itinerary $\chi(K_j)$. 

Since $g=g^t$ is the full family, there is a parameter interval $T=[t_0,t_1]$
such that the critical value $g c$ runs though the interval $L_0$ from
one boundary point to another as $t$ runs through $T$.  
Since $\chi$ is unimodal (and $\Gamma_g$ is positively oriented),
 $L_0$ lies on the left of all other intervals $L_j$. 
Hence for $t\in T$ we have $g I_0\supset L_i,\; i\neq 0$. 

Let us now consider the pull-back $K_j\subset I_0$ of  the intervals $L_j$ 
 by the central branch in such a way that $K_j$ lies on the right of $c$
if and only if $j>0$. Then the first return map $h$ on $\cup K_j$ has
type $\chi$. Moreover, as $t$ runs through $T$, the critical value $h c$
runs all way through $I^0$, so that we have a full family. \QED

Let $g=g_1\in {\cal T}$. Assume that we can subsequently construct a sequence
of return maps of class ${\cal T}$
$$g_n : \cup I^n_k\rightarrow I^{n-1}_0$$  where
 $I^0_0\supset I^1_0\supset....$ is a sequence of the first kids.
Let  $\chi_n,\; n=1,2,...$, be the corresponding
sequence of return types. The following statement says
that these types can be combined independently. (``Do whatever you want".)

\proclaim Lemma 2.7. Let $g^t\in {\cal T}(q,p,\epsilon),\; 0\leq t \leq 1,$ 
be a full one-parameter family, $\Gamma_0\equiv \Gamma_g$.
Let $\chi_n: \Gamma_n\rightarrow \Gamma_{n-1},\; n=1,2,...$ be any
sequence of  admissible homomorphisms.
Then there is a map $g=g^t$ which admits a sequence of return maps
$g_n$ of type $\chi_n$.

\noindent {\bf Proof.} By the previous lemma  we can subsequently construct 
a nested  sequence  of parameter intervals $T_1\supset T_2\supset...$ in 
such a way that the return types to   $T_k$ are required up to the level $n$,
while the retun maps $g_{n+1}^t,\; t\in T_n$,
 to the $n$th level form a full family.
Intersecting these intervals we get a required parameter value. \QED

 In particular, we can fix the intervals  $I_k$ and consider the 
``standard" family
of maps which are linear on the non-central intervals (and hence don't
depend on parameter) and quadratic on the central interval. Then by
changing this central branch we can obtain a map of any type
$\{\chi_n\}$.

\medskip\noindent {\bf The return graph.} Let us consider a sequence
$$\matrix{...&\rightarrow&\Gamma_2&\rightarrow &\Gamma_1&\rightarrow&\Gamma_0\cr
             &    \chi_3 &        &    \chi_2  &        &   \chi_1   &  \cr } $$
 of admissible homomorphisms of spin semigroups.
 We can associate to it the following
graded graph $\Delta$. Put on the $n$th level of $\Delta$ the generators
$\{ I^n_k\}$ of $\Gamma_n$. Connect a vertex $I^n_k$ with a vertex 
$I^{n-1}_j$ of
the previous level by $l$ edges if $I^{n-1}_j$ has multiplicity $l$
in the word $\chi_n(I^n_k)$. This graph contains the full ``abelian" 
information about the sequence of homomorphisms.

Observe that each vertex of level $n$ is connected by a simple edge with
the central vertex $I^{n-1}_0$ of the previous level. Let us call the graph
(and the sequence $\{\chi_n\}$)
{\sl irreducible} if for any vertex $I^n_k$ there is a path down to a
central vertex $I^{n+s}_0$
(such that we go strictly downstairs along this path).

Let us consider now a sequence of return maps $g_n$. Its type is given
by a sequence of homomorphisms $\chi_n$ which provides us with a graph $\Delta$.
A vertex $I^n_k$ is connected with $I^{n-1}_j$ by $l$ edges 
if the $g_{n-1}$-orbit of
$I^n_k$ passes through $I^{n-1}_j$ $l$ times before the first return back to the
central interval $I^{n-1}_0$. (Such a graph was introduced by Marco Martens.)
The following statement is immediate:

\proclaim Lemma 2.8. The graph $\Delta$ is irreducible if and only if the
critical orb($c$) crosses all intervals $I^n_k$. \QED

\medskip\noindent{\bf Realization of all types.}
 Let $g\in {\cal T}$ have  a recurrent critical
point such that $f|\omega(c)$ is minimal. Assume also that $g$ does not
admit a unimodal renormalization. Denote this class of maps by 
${\cal T}_{\min}$. Let us consider the nested
sequence of the first grandkids $I^0\supset I^1\supset...$. Let us define
a  sequence of return maps
$$g_n : \cup I^n_k\rightarrow I^{n-1},\quad I^n\equiv I^n_0$$
in such a way that  we select only those intervals $I^n_k$ which cover
the critical set. The number of these intervals is finite since $\omega(c)$
is minimal. We call $g_n$ a {\sl pre-renormalization}  of $g_{n-1}$.
So, we obtain a sequence $\bar\chi_g=\{\chi_n\}$ of 
admissible homomorphisms with an irreducible graph $\Delta_g$. 
Let us call this sequence a {\sl type} of $g$.  

\proclaim Lemma 2.9.  Let $g^t\in {\cal T}(q,p,\epsilon),\; 0\leq t \leq 1$ 
be a full one-parameter family, $\Gamma_0\equiv \Gamma_g$.
Let $\bar\chi=\{\chi_n : \Gamma_n\rightarrow \Gamma_{n-1}\}_{n=1}^\infty$ 
be any irreducible sequence of  admissible homomorphisms 
such that $\chi_n (I^n)\neq (I^{n-1})$ for infinitely many $n$.
Then there is a map $g^t\in {\cal T}_{\min}$  of type $\bar\chi$. 

\noindent{\bf Proof.} By Lemma 2.7 we know that there is a map $g^t$ which
admits a sequence of return maps of type $\chi_n$. Since 
$\bar\chi$  is irreducible, this is a sequence of pre-renormalized maps
(by Lemma 2.8). Hence $\bar\chi_g=\bar\chi$. 

 This map does not admit a unimodal renormalization
since $\chi_n (I^n)\neq (I^{n-1})$ for infinitely many $n$.
 By Lemma 2.3 diam$(I^n)\to 0$. Since the
critical point returns to all central intervals $I^n$, it is recurrent. 
In order to see the minimality property, let us consider the union of
orb($I^n_k$) until their first return back to $I^{n-1}$. This provides
us with a finite covering of $\omega(c)$ by the pull-backs of $I^{n-1}$.
It follows that the orbit of any point $x\in\omega(c)$ crosses $I^{n-1}$
and hence $\omega(x)=\omega(c)$. \QED


\medskip\noindent {\bf Combinatorial model of the critical set.}
 Let $\Gamma_*^m$ be a sub-semigroup
of $\Gamma^m$ generated by non-central intervals $I_k,\; k\neq 0$.
Consider  a space $\Omega$ of all finite and infinite sequences
$\bar\delta=(\delta_1,\delta_2,...)$
such that $\delta_n\in \Gamma^{l(n)}_*$ with
 a {\sl strictly increasing} sequence of levels $l(n)$.
 Let us define a map $\sigma: \Omega\rightarrow \Omega$ in the following
fashion.
Set $n=l(1)$.
If the word $\delta_1$ has length greater than 1 then just forget 
the first symbol of this word. 
If $\delta_1=(I^n_k)$  and
$\chi(\delta_1)=(I^{n-1}_0)$ then forget $\delta_1$. Finally, if $\delta_1$
has only one symbol and $\chi(\delta_1)\not= (I^{n-1}_0)$ then
replace $\delta_1$ with $\chi_*(\delta_1)$ where
 $\chi_*(\delta_1)$ coincides with
$\chi(\delta_1)$ except for the last symbol $I_0^{n-1}$ which is dropped.

Given a map $g\in {\cal T}_{\min}$, we can  associate to a point $x\in\omega(c)$
an element of the space $\Omega$ in the following way.
Let $x\in I^{n-1}\ssm I^n$. Then consider the itinerary $\delta_1$ of $x$ 
 through the intervals of level $n=l(1)$ until the first moment $t$
when it lands in the central interval $I^n$.
Then find a level $l(2)=m>n=l(1)$ 
such that $y=g_n^{\circ t}x\in I^m\ssm I^{m-1}$. Coding now $y$ in a
 similar way we will find $\delta_2$, etc.
Since ${\rm diam} (I^n_k)\to 0$ as $n\to\infty$, this provides us with a 
homeomorphism between $\omega(c)$ and $\Omega$ conjugating $f$ and $\sigma$.
So, we have the  following lemma. 

\proclaim Lemma 2.10. Two maps $f,g\in {\cal T}_{\min}$ have the same 
type, $\bar\chi_f=\bar\chi_g$, if and only if their restrictions on the
critical sets are topologically
conjugate  by an orientation preserving homeomorphism  respecting the 
critical points.  \QED

\medskip\noindent {\bf Pull-back argument.}

\proclaim Lemma 2.11. Two maps $f$ and $g$ 
of class ${\cal T}_{\min}$  are topologically
conjugate by an orientation preserving homeomorphism 
if and only they have the same types, $\bar\chi_f=\bar\chi_g$. 

\noindent{\bf Proof.} Let us set $g\equiv {\tilde f}$ and mark all objects
related to $g$ with the tilde. By Lemma 2.9, there is an orientation
preserving conjugacy $h:\omega(c)\rightarrow\omega(\tilde c)$ putting
$c$ to ${\tilde c}$. Let us continue it to an orientation preserving
homeomorphism $h: (J, \cup I_k)\rightarrow ({\tilde J},\cup{\tilde I}_k)$ 
of the base intervals.

Consider now a nested sequence of inverse images 
 $J^{(k)}=f^{-k} J,\; k=0,1,...$
shrinking down to the Julia set $K(f)$, $J^{(1)}\equiv \cup I_k$ 
(and define ${\tilde J}^{(k)}$ similarly).
Since $h(c)=\tilde c$, we can lift $h$ via the folding maps $f$ and $\tilde f$
to an orientation preserving homeomorphism:

$$\matrix{      &     h_1     &             \cr
           \quad J^{(1)} & \rightarrow &\tilde J^{(1)}  \cr
          f\downarrow&       &\quad\downarrow \tilde f\cr
          \quad J        &  \rightarrow &\tilde J\cr
                   &     h        &        \cr}$$ 

Then $h_1|\partial J^{(1)}=h$, and hence $h$ provides a continuation of $h_1$
 to the whole base interval $J$. 

Let us show that $h_1|\omega(c)=h|\omega(c)$. Indeed, since $h$ is a 
conjugacy on $\omega(c)$, the points $h(x)$ and $h_1(x)$ either
coincide or are $\tilde c$-symmetric. Since both of maps are orientation 
preserving, they must coincide.

Now let us repeat the construction and pull $h_1$ back, etc. In such a way
we will construct a sequence $h_n: J\rightarrow\tilde J$ of orientation
preserving homeomorphisms such that

\noindent (i) ${\tilde f}\circ h_n=h_{n-1}\circ f$;

\noindent (ii) $h_n$ agrees with $h$ on the critical set;

\noindent (iii) $h_n|J\ssm J^{(n)}=h_{n-1}$.

Hence there is a pointwise limit $H=\lim h_n$ which is an orientation
preserving homeomorphism $J\ssm K(f)\rightarrow {\tilde J}\ssm K({\tilde f})$
conjugating $f$ and ${\tilde f}$. Since the  both Julia sets
 $K(f)$ and $K({\tilde f})$ are Cantor,  $H$ can be automatically
continued across these sets.  \QED

\medskip\noindent {\bf Example: the Fibonacci recurrence.}
In this case we have two intervals $I^n_0$ and $I^n_j,\; j\in\{-1,+1\},$ 
on all levels.
The returns are high on all levels, and 
 $\chi(I^n_0)=(I_j^{n-1},I_0^{n-1})$, $\chi(I_j^n) =(I_0^{n-1})$. 
This data and the spin structure on the first level
(four possibilities) uniquely determine the admissible spin structures on the
further levels and the choice between $j=1$ or $-1$ for the 
non-central intervals
$I^n_j$ Namely, one can see that on two subsequent levels the $I^n_j$ lie 
on one side of $c$, and then on the next two  levels they lie on the opposite
side, etc (the first level can be the only exception).

 We see that there are four Fibonacci types. But actually two pairs of them 
are  the same up to the choice
of the orientation of the base interval. So, only two types are left
which differ by the spin of the non-central branch.  The renormalization
interchanges these types.

Let us remark in conclusion that the name ``Fibonacci" comes from the
observation that the first return time $t_n$ of the critical point to
the $n$th level satisfies the Fibonacci recurrent equation
$t_{n+1}=t_n+t_{n-1}$.

\medskip\noindent{\bf Renormalization.} In the subsection ``Realization
of all types" we defined  the pre-renormalization $g_1$ of a map
 $g\in{\cal T}_{\min}$.
The {\sl renormalization} $Rg$ is obtained from $g_1$ just by 
rescaling of the base interval $I^n_0$ (for $g_1$) to the original interval
$J$. If we repeat this procedure, we will see that the $n$-fold 
renormalization $R^{\circ n} g$ is obtained from $g_n$ by rescaling $I^n_0$.
The type of the renormalization is defined as  the type of the corresponding
pre-renormalization. Now we can summarize the above results as follows: 

\proclaim Theorem 2.12. Let $g\in {\cal T}$.  Then $g$ is infinitely
renormalizable with $R^{\circ n} g\in {\cal T}_{\min}$. The topological type
of $g$ is uniquely determined by the sequence $\bar\chi_g$ 
of renormalization types.
Any sequence of admissible types 
$$\matrix{...&\rightarrow&\Gamma_2&\rightarrow&\Gamma_1&\rightarrow &\Gamma_0\cr
             &   \chi_3   &        &  \chi_2   &        &  \chi_1    &  \cr}$$
can be realised in any full one-parameter family $g^t$ with $\Gamma_g=\Gamma_0$.

\medskip\noindent {\bf Maps of infinite type.} In the case when the 
critical point is recurrent but $\omega(c)$ is not minimal, we still
can renormalize the map on a nice interval
by taking the first return map, keeping only the
intervals intersecting $\omega(c)$ and rescaling. However, now the 
domain of the $n$-fold renormalization
consists of infinitely many intervals starting from some level.
Still one can develop a similar combinatorial theory but we don't need it.
Let us only agree that $I^n\equiv I^n_0$ is still the central interval,
and the intervals $I^n_k$ with $k>0/k<0$ lie on the right/left of $c$.


\medskip\noindent{\bf Getting started.} Let $f: L\rightarrow K,\; L\subset K,$
 be a unimodal map, exact on $[fc, f^{\circ 2}c]$. Take a nice interval
$J\equiv I^0\subset    [fc, f^{\circ 2}c]$,
 and renormalize $f$ on $J$. We obtain a
map $g\in {\cal T}^{\infty}$ which does not admit unimodal renormalizations.
This is our object to study.

\bigskip
 \centerline{\bf \S 3. Estimates of Poincar\'{e} lengths and 
 scaling factors.}\medskip

\noindent{\bf Hyperbolic line and asymmetric Poincar\'{e} length.}
Let us consider an interval
$L=[a, b]$ as a {\sl hyperbolic line} with the Poincar\'{e} metric
$2 dx/(x-a)(b-x)$. 
Let  $ G$ be a subinterval of $L$, and $U$ and $V$ be the
components of $L\ssm G$ (see Figure 1).  We will use the following notations
for the Poincar\'{e} length of $G$ in $L$:
$$P(G)\equiv P(G | L) \equiv P(U,V)=\log (1+{|G|\over |U|}) +
\log (1+{|G|\over |V|}).$$ 
Note that the bigger the space is around
$G$ in $L$, the smaller Poincar\'{e} length $P(G|L)$ is.

\medskip
\centerline{\psfig{figure=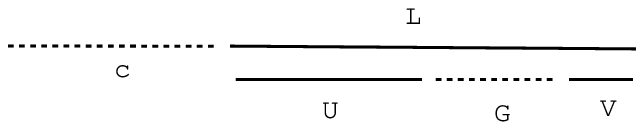,height=.8in,width=2.5in}}
\centerline{Figure 1}\medskip

Let us now state the main analytic tools of real one-dimensional dynamics
(see [MS]).
They are  called ``the Schwarz Lemma"  and ``the Koebe Principle"
by analogy with the classical facts in 
geometric function theory.

\proclaim Schwarz Lemma. Any diffeomorphism $h: L\rightarrow L'$
with positive Schwarzian derivative contracts the Poincar\'{e} metric.
\QED

Hence, given an interval $G\subset L$ and its image $G'= hG$ we have:
$P(G | L)\leq P( G' | L')$.  So, if we have a
definite space in $L$ around $G$ then we have also a definite space in
$L'$ around $G'$.

Given a diffeomorphism $h: I\rightarrow I'$, let us call
$$\max_{x,y\in I}\log{|f'(x)|\over |f'(y)|}$$
its {\sl distortion} or {\sl non-linearity}. The case of zero non-linearity
corresponds to linear maps.

\proclaim Koebe Principle. Let  $h: (L,I)\rightarrow (L',I')$ 
be a diffeomorphism with positive Schwarzian derivative, $r=P(I|L)$. 
Then the non-linearity of $h$ on $I$ is bounded by a constant $C(r)$
independent of $h$. Moreover, $C(r)=O(r)$ as $r\to 0$. \QED

\noindent
\proclaim Lemma 3.1. \it If $L\supset T\supset I$ then
$ P(I | L) \leq  {1\over 2} P(I | T)\, P( T | L).$

\noindent{\bf Proof.} Since the Poincar\'{e} metric is invariant under
M\"{o}bius transformations, we can normalize the intervals in the 
following way:  $L=[-1,1]$, $T=[-\lambda,\lambda],\; 0<\lambda<1$.
Let us consider the map $h: x\mapsto\lambda x$. The calculation
shows that $\|D h(x)\|_L\leq\lambda$ in the $L$-Poincar\'{e} metric 
(with equality at 0). On the other hand, $h: L\mapsto T$ is an isometry
from the $L$-Poincar\'{e} metric to the $T$-one. Hence
$P(I|L)\leq\lambda\, P(h^{-1}I|L)=\lambda\, P(I|T).$
Observe finally that $2\lambda\leq P(T|L)$.  \QED


\smallskip\noindent {\bf Remark.} We will actually use Lemma 3.1 in a slightly
  different  form:
$$ P(I | L) \leq  P(I | T)\, P^*( T | L)$$
where  $P^*(T|L)=\min \{P(T|L), 1\}$.

Suppose we have a map $g$ with a single non-degenerate critical point $c$.
If the interval int$G$  does not contain $c$,
 then let us introduce the {\it asymmetric Poincar\'{e}
length}  defined as

$$Q(G)\equiv Q(G | L) \equiv Q(U,V)= \log (1+{|G|\over |U|})+
   {1\over 2} \log (1+{|G|\over |V|}),$$
provided $U$ is closer to $c$ than $V$. Clearly, $Q(G)<P(G)$. The coefficient
1/2 is related to the exponent 2 of the critical point $c$. It turns out
that the asymmetric Poincar\'{e} length behaves more regularly under
renormalizations of a quasi-quadratic map than the usual Poincar\'{e} length.  

\medskip\noindent {\bf Parameters.}
From now on we will assume without change of notations that all maps
{\sl of class ${\cal T}$ have negative Schwarzian derivative and a 
non-degenerate critical point} $c$. Let us consider a map 
$g\equiv q_1\in {\cal T}$  which {\sl does not
admit unimodal renormalizations, and has a recurrent critical point}.
Then we can construct the sequence of pre-renormalized maps
$$ g_n: 
\bigcup I^n_k\rightarrow I^{n-1}$$
with the central intervals $I^n_0\equiv I^n$ shrinking down to $c$. 
In this section we don't assume that
$\omega(c)$ is minimal, so that we allow infinitely many intervals $I^n_k$
on all sufficiently high levels $n$.

By the {\sl gap} between intervals $U$ and $V$ 
we mean the bounded connected component of ${\bf R}\ssm (U\cup V)$.
Let us introduce the following parameters: 

\smallskip\noindent - $K_n$ is the infimum of asymmetric Poincar\'{e}
lengths $Q(I^n_s, I^n_t)$ of the gaps between intervals 
of level $n$ with  $st\geq 0$; 
 
\smallskip\noindent -  $\mu_n=|I^{n}|/ |I^{n-1}|$ is 
the scaling factor on level $n$;

\smallskip\noindent -
 $\lambda_n=\max_{i\neq 0}\, [I^{n-1} : I^n_i]$ is
the maximal Poincar\'{e} length
of the non-central intervals of \break $\quad$ level $n$. 
Set $\lambda_n^*=\min(\lambda_n,1).$

 Further, let
$$\alpha_n=\sup_{k\neq 0}{I_k^n\over {\rm dist}(I_k^n,c)} $$
be a parameter
which controls the non-linearity of the quadratic map $\phi(x)=(x-c)^2$ on 
non-central intervals.

By means of a $C^3$-small change of variable (near the critical value)
 we can make $f$ 
{\sl purely quadratic} in a neighborhood of the critical point $c$.
Then  $g_n$ can be decomposed  in the following way:
 $$g_n|I^n_i=h_{n,i}\circ\phi,\eqno (3-0)$$
where  $\phi(x)=(x-c)^2$ is the quadratic map and
$h_{n,i}$ is a diffeomorphism with negative Schwarzian derivative 
of an appropriate interval onto $J$.
Let us consider one more parameter:

- $\rho_n$ is the maximal distortion
of $h_{n,i}$, and call it the {\sl distortion parameter}.

Let us remember that  $\kappa: {\bf N}\rightarrow {\bf N}$ numbers the
cascades of central returns (see \S 2).
The goal of this section is to prove the following: 

\proclaim Conditional version of Theorem B. 
There exist $\bar K$ and $\bar \mu>0$
(independent of a map) with the following property.
 If on some level $N$, $K_N>\bar K$ or $\mu_N<\bar \mu$
   then  $K_n\to\infty$. Moreover,
there exist positive constants $A$, $C$ and $\sigma$ such that
$$K_n\geq A \,\kappa(n), $$ 
$$\mu_{n+1}\leq C\exp(-\sigma\kappa(n)), \; n\in {\cal L},$$
$$\rho_n\leq C\exp(-\sigma\kappa(n)).$$

\noindent {\bf Remark.}  Observe that the  scaling factors $\mu_n$
are exponentially small only
 on the special subsequence of levels (outside the long cascades of
central returns). 
 However, the estimates of the Poincar\'{e} lengths of
the gaps, as well as the non-linearity control hold on all levels.\smallskip 

Let us fix constants
$\bar\mu$ and  $\bar\rho$  and $\bar K>0$ such that $1>\bar\mu>\bar\rho>0$.
We will assume (until the subsection ``Cascades of central returns")
that the following estimates hold:

$$ \qquad\rho_n<\bar\rho,
 \qquad \mu_n<\bar\mu,\qquad K_n>\bar K.  \eqno (3-1) $$

So, $\bar\rho$ controls the distortion, 
 $\bar\mu$ controls the scaling factors, and $\bar K$ controls the Poincar\'{e}
lengths of the gaps.
In what follows all constants depend on $\bar\mu$ (and actually on $\bar\rho$
which becomes non-important because we keep $\bar\mu>\bar\rho$) but not
on the particular map.  
 Sometimes we will abuse notations  using the
same letter for different constants.

 Let us also fix small constants $\bar\lambda$ and $\bar\alpha$ which
separate range of small values of parameters $\lambda_n$ and $\alpha_n$
from big ones. 

\medskip\noindent {\bf Strategy.}
Our strategy is the following. Let us consider two intervals 
$U'=I^{n+1}_k$ and $V'=I^{n+1}_j$ such that the gap $G'$ between them does
not contain the critical point. 
Let us push these intervals forward by $g_n$ until the first moment $p$ when 
$U= g_n^p U'\subset I^n_t$ and $V=g_n^p V'\subset I^n_s$ lie in different 
intervals of level $n$, that is $s\neq t$. Then loosely speaking, the 
Poincar\'{e} length of the gap between $U$ and $V$ can
be estimated from below by $2 K_n+\chi$ with an absolute constant $\chi>0$.
Pulling this back by an almost quadratic map we get an estimate of the
asymmetric Poincar\'{e} length of the gap between $U'$ and $V'$, namely
$Q(U,V)\geq K_n+\chi/2$.

The argument  depends on the positions of the intervals $I^n_t$ and $I^n_s$.
The Fibonacci-like situation when one
of these intervals is central is the main one to look at
(see Lemmas 3.3 and 3.9).
In all other cases the estimates are actually getting better.

\medskip\noindent {\bf Estimates of $P(U,V)$.}
This will occupy lemmas 3.2 through 3.7.

Let us fix a level $n$ and temporarily drop the index  $n$ in all notations
 so that $I^n\equiv I,\; g_n\equiv g,\; \mu_n\equiv \mu$ etc.
However, let $I^{n-1}\equiv J$.
Let us take a non-central interval $I_t, \; t\neq 0$, of level $n$, 
and consider an interval 
$U\subset I_t$ such that $g^l U=I,\; g^i U\cap I=\emptyset,\; i=0,1,...,l-1.$
Sometimes we will write $l=l_U$. 
Let $L$ and $R$ be the components of $I_t\ssm U$ with $L$ closer to $c$ than
$R$ (see the following figure).

\centerline{\psfig{figure=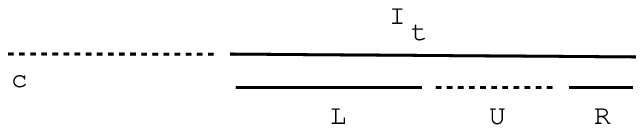,height=2in,width=2.5in}}
\centerline{Figure 2}\medskip

Let $J^+$ and $J^-$ be the components of $J\ssm I$.

\proclaim Lemma 3.2.  The following estimates hold:
$$ {|U|\over |L|}+{|U|\over |R|}\leq 4\mu(1+O(\mu)). \eqno (3-2)$$
If $l=1$ then
$$ {|U|\over |L|}\leq 2\mu (1+O(\rho+\mu))\eqno (3-3) $$ and
$$  C(\rho,\mu)^{-1}\leq {|R|\over |L|}\leq C(\rho,\mu), \eqno (3-4)$$
with $C(\rho,\mu)= \sqrt{2}(1+O(\rho+\mu))$.
If $l>1$ and then
$$  {|U|\over |L|}\leq 4\mu\lambda^*(1+O(\mu))\qquad  and\qquad  
{|U|\over |R|}\leq 4\mu\lambda^*(1+O(\mu)). \eqno (3-5)$$


\noindent{\bf Proof.}  There is an interval $T$ such that 
$U\subset T\subset I_t$ which is diffeomorphically mapped by $g^l$ onto $J$.
By the Schwarz lemma,
$$P(U|I_t)\leq P(U|T)\leq P(I|J)=4\mu (1+O(\mu)). \eqno (3-6)$$
But $$P(U|I_t)=({|U|\over |L|}+{|U|\over |R|})\,(1+O(P(U|I_t))), \eqno (3-7)$$ 
provided  there is an a priori bound on $[I_t:U]$. 
The last two estimates imply (3-2) .

In order to get (3-3) let us make use of the decomposition (3-0):
$${|U|\over |L|} < {|\phi U|\over |\phi L|}<
{|I_0|\over  |J^+|} (1+\rho)=2\mu(1+O(\rho+\mu)).$$

Estimate (3-4) follows from the fact that $g$ is the composition of a
quasi-symmetric map $\phi$ and a diffeomorphism with distortion $\rho$
(the $\sqrt{2}$ comes as the qs norm of $\phi^{-1}$).

Suppose now that $l>1$. 
Then $g^{l-1}$ maps the interval $T$ introduced above onto a
non-central interval $I_s,\; s\neq 0$. Hence there is another interval 
$T'$ in between $T$ and $I_t$ which is mapped by $g^{l-1}$ onto $J$. Hence
$P(T|I_t)\leq P(T|T')\leq P(I_s|J)\leq\lambda.$
By Lemma 3.1 and estimate (3-6),
$$P(U|I_t)\leq P(T|I_t) P(U|T)\leq  4\mu\min(\lambda,1) (1+O(\mu)),\eqno (3-8)$$
and the estimates (3-5) follow from (3-7) and (3-8). \QED \medskip

We will use the sign $\prec$ or $\succ$ if an estimate holds up to 
$O(\mu+\rho)$,
and a sign $\approx$ if an equality holds up to $O(\mu+\rho)$ (provided
$\mu\leq\bar\mu, \rho\leq\bar\rho$). 

Let $U\subset I_t$ be as above, $H$ be the gap between $I_0$ and $I_t$,
$G$ be the gap between $I_0$ and $U$ (see Figure 3). 
Let $P(H)$ denote  the Poincar\'{e}
length of $H$ in $I_0\cup H\cup I_t$, and $P(G)$ denote the Poincar\'{e}
length of $G$ in $I_0\cup G\cup U$. Notations $Q(H)$ and $Q(G)$ mean
the asymmetric Poincar\'{e} lengths of the same pairs of intervals. 
Let $J^+$ be the component of $J\ssm I$ containing $I_t$.

\bigskip
\centerline{\psfig{figure=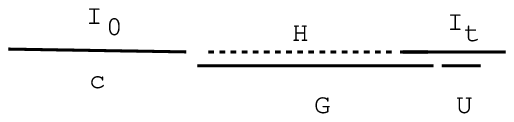,height=.6in,width=2.5in}}
\medskip\centerline{Figure 3}\medskip

\proclaim Lemma 3.3. 
If $l=1$ then there is an absolute constant $\chi>0$ such that
$$P(G)\succ 2 Q(H)+\chi. \eqno (3-10)$$
If $l>1$ then
$$P(G)\succ 2 Q(H)+\log^+{1\over\lambda}-\log 2. \eqno (3-11)$$

\noindent
{\bf Proof.} We have
$$\eqalign{P(G)\geq\log(1+{|H|\over |I|}) + \log(1+{|H|+|L|\over |U|})=\cr
=\log(1+{|H|\over |I|}) + \log(1+{|H|\over |L|+|U|})+
\log(1+{|L|\over |U|}).\cr}  \eqno (3-12)$$
The middle term is evidently bounded from below by  $\log(1+|H|/ |I_t|)$.
As to the last term, then by (3-2) we have:
$$\eqalign{\log(1+{|L|\over |U|})\succ \log{1\over 4\mu}
\approx\log(1+{|J^+|\over |I|})-\log 2\geq\cr
\geq\log(1+{|H|\over |I|})+\log{|J^+|+|I|\over |H|+|I|}-\log 2,\cr}
 \eqno (3-13)$$
and the estimate $P(G)\succ 2 Q(H)-\log 2$ follows.

Let now $l=1$. Then we can use (3-3) instead of (3-2), and $-\log 2$
disappears in the last estimate. We can also improve the estimate of the 
middle term of (3-12) as follows. Because of (3-3) and (3-4),
there is a $\tau >1$ such that
          $|I_t|\geq \tau (|L|+|U|). $
Hence
$$ \log(1+{|H|\over |L|+|U|})\succ
\log(1+{|H|\over |I_t|})+\chi, \eqno (3-14)$$
provided $H$ is not tiny as compared with $I_t$.
On the other hand if $H$ is tiny as compared with $I_t$ then 
$J^+$ is big as compared with $H$ and $H$ is big as compared with $I$
(namely, $\log(|H|/|I|)\geq\bar K-$[{\sl a tiny term}]). Hence the second term
$$\log{ |J^+|+|I|\over |H|+|I|}$$ in (3-13) is big, and suppresses $-\log 2$.
These yield (3-10).

Finally, if $l>1$ then we can improve (3-13) by using (3-5) instead of
(3-2).  \QED


Now together with the above pair of intervals $U\subset I_t$, let us consider a 
similar pair  $V\subset I_s, \; s\neq 0,$ with $V$ to be a monotone pull-back
of $I$. Let us assume that both pairs lie on the same side of $c$, and the
latter one is closer to $c$ than the former (see Figure 4).
Let $G$ be the gap between $V$ and $U$, $H$ be the gap between $I_s$ and
$I_t$. 

\centerline{\psfig{figure=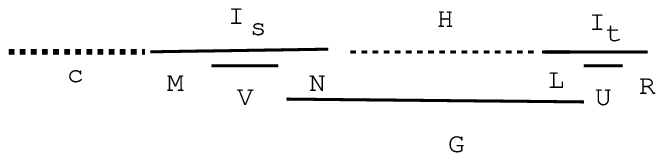,height=2.1in,width=3in}}
\centerline{Figure 4}\medskip


\proclaim Lemma 3.4.
If $I_s$ and $I_t$ are non-central intervals lying on the same side of $c$ then 
$$P(G)\geq P(H)+2\log{1\over \mu}-O(1). \eqno (3-16)$$

\noindent
{\bf Proof.} Let $L$ and $R$ be  the components of $I_t\ssm U$ 
as defined above, while
$M$ and $N$ are the components of $I_s\ssm V$. Then we have:
$$\eqalign{P(G)=\log(1+{|H|+|N|\over |V|}) +\log(1+{|H|+|L|\over |U|})=\cr
=\log(1+{|H|\over |N|+|V|})+\log(1+{|N|\over |V|})+\cr
+ \log(1+{|H|\over |L|+|U|})+\log(1+{|L|\over |U|}).\cr}\eqno (3-18)$$
The sum of the first and the third terms of (3-18) is certainly greater than
$P(H)$. Estimating the second and the last terms by (3-2), we get (3-16). \QED



The following lemma will allow us to handle the case when the non-linearity
of $\phi$ is not small.

\proclaim Lemma 3.5. 
The following estimate holds:
$\mu=O((1+{1\over\alpha}) e^{- K}$.

\noindent
{\bf Proof.} Let us select an interval $I_k$ for which
$\alpha=|I_k| / {\rm dist}(I_k,c).$
Let $W$ be the gap between $I$ and $I_k$. Then we have:
$$K\leq P(I,I_k)\prec\log{1\over 2\mu}+\log(1+{1\over\alpha}),$$
and the conclusion follows.  \QED

We will need the following lemma to analyze cascades of central returns.

\proclaim Lemma 3.6. Under the circumstances of Lemma 3.4,
if $H\geq {\rm dist}(H,c)$ then\break
$P(G)\geq (5/2) K-O(1).$

\noindent
{\bf Proof.} Let $W$ be the gap between
$I$ and $I_s$ and $X$ be the gap between $I$ and $I_t$.
 Let us start with formula (3-18). Because of the  assumption of the lemma, 
we can estimate from below  half of its first term by
$(1/2)\log (1+|W|/|I_s|)$, and  half of the third term by
$(1/2)\log(1+|X|/|I_t|)-\log 2$. 
The sum of the other halves we estimate as $(1/2)P(H)$.

Remember that $J^+$ denotes a component of $J\ssm I$.
The second and the last terms we $\succ$ estimate by (3-2) as
$\log(1+|I|/|J^+|)$ which is greater than both
$ \log(1+|I|/|W|)$ and $\log (1+|I|/|X|)$.
Taking all these together, we get
$$P(G)\geq Q(I,I_s)+Q(I,I_t)+{1\over 2} P(H)-O(1)\geq (5/2)K-O(1).$$ \QED
Finally, let us consider the case when $I_s$ and $I_t$ lie on the opposite
sides of $c$.  

\proclaim Lemma 3.7. If $I_s$ and $I_t$ lie on the opposite sides of $c$
then $P(G)\succ (5/2)K.$

\noindent
{\bf Proof.} The argument is the same as in the previous lemma. The point is
that now we automatically have $|H|\geq |W|$ and $|H|\geq |X|$ where 
as above $W$ denotes
the gap between $I$ and $I_s$, and $X$ denotes the gap between $I$ and $I_t$.
\QED

\medskip\noindent {\bf Quadratic pull-backs.} 
Let us start with a lemma which says that the square root map divides the
Poincar\'{e} length at most by 2.

\proclaim Lemma 3.8. Let us consider a quadratic map $\phi: x\mapsto (x-c)^2$.
Let $U$ and $V$ be two disjoint intervals lying on the same side of $c$,
$V$ being closer to $c$ than $U$. Then 
$$Q(V,U)>{1\over 2} P(\phi V,\phi U).$$

\noindent
{\bf Proof.} We can assume that $c=0$ and $V$, $U$ lie on the right of $c$.
 Let $V=[v,a]$, $U=[b,u]$. Then 
$$P(V,U)-{1\over 2} P(\phi V,\phi U)=$$
$$={1\over 2}\left( \log{a+v\over a-v} -\log{b+v\over b-v}\right)
 +{1\over 2}\left( \log {u-a\over u-b} + \log {u+b\over u+a}\right)>
  {1\over 2}\log{u-a\over u-b},$$
which is exactly what is claimed.  \QED 

 In Lemmas 3.3 - 3.7
we have estimated  the 
Poincar\'{e} length of the gap $G$ between $U$ and $V$.
Now we are going to use Lemma 3.8  in order to estimate
the asymmetric Poincar\'{e} length of the gap between $U'$ and $V'$.
Again let us start with the situation when one of the intervals, say $I_s$,
is central (as in Lemma 3.3). As above, $H$ denotes the gap between $I_s$
and $I_t$. Set $T'= U'\cup G'\cup V',\;  T=g^p T'=U\cup G\cup V.$

\proclaim Lemma 3.9.   
 Under the circumstances just described
there is a constant $\chi>0$ such that 
$$Q(G')\succ Q(H)+\chi, \eqno (3-22)$$ or
$$Q(G')\succ K+{1\over 2}\log^+{1\over \lambda}-O(1). \eqno (3-23)$$

\noindent
{\bf Proof.} 
{\sl Case 1.}  Let $p=1$.

\noindent  Let us use representation (3-0) of $g|I$ as the quadratic map
$\phi$ postcomposed by a diffeomorphism $h$ with distortion $\rho$.
 Pulling $G$ back by $h$ and then by $\phi$ (making use of Lemma 3.8), we see
that  $Q(G')\succ (1/2) P(G)$. Together with Lemma 3.3
 this yields the claim.  
\smallskip

{\sl Case 2.} Let $p>1$.

\noindent Then let us consider the intervals $\tilde U = g^{p-1} U',\;
\tilde V =g^{p-1} V' $ and $\tilde G =g^{p-1} G' $. All three of them
belong to the same interval of level $n$, say $I_j$. 

 Because of (3-0) we can consider the following decomposition:

$$ \psi\circ \phi|T', \eqno (3-25) $$
where $\psi$ is a diffeomorphism onto $J$.
Remember that the non-linearity of the
quadratic map $\phi|I_j$ is controlled by the quantifier $\alpha$.
Let us take a small $\bar\alpha>0$ and consider  several subcases.\smallskip

{\sl Subcase (i).} Assume $\alpha< \bar\alpha$.

\noindent This implies that
 $g|I_k,\; {k\neq 0}, $ is an
expanding map with the expansion  $>2$ and small non-linearity.
Then the diffeomorphism $\psi$ in (3-25) has small non-linearity as well.
Together with Lemmas 3.3 and 3.8 this yields the desired estimates.
\smallskip

Let $S$ be the gap between $I$ and $U$.

{\sl Subcase (ii).} Assume that $|S|\leq \bar\alpha |J|$.

\noindent Then $P(T|J)=O(\bar\alpha+\mu)$ (make use of Lemma 3.2). 
Hence $\psi|\phi T'$ has small non-linearity,
 and the result follows.\smallskip

\smallskip
{\sl Subcase (iii).} Finally, assume that 
$$\alpha> \bar\alpha\eqno (3-26)$$
      and $${|S|\over |J|} > \bar\alpha . \eqno (3-27)$$
Then Lemma 3.5 and (3-26) imply
$${|I|\over |J|}=\mu=O(e^{-K}). \eqno (3-28)$$
Together with (3-27) this implies
$${|I|\over |S|}=O(e^{-K}). \eqno (3-29)$$
Let $J^-$ be the component of $J\ssm I$ disjoint from $S$. Given an interval
$X\subset I$, let $\tilde X$ denote its pull-back  by $\psi$. 
Pulling  the interval $J^-\cup I\cup S$ back by $\psi$,  we get by
(3-28), (3-29) and the Schwarz lemma 
$${|\tilde I|\over |\tilde S|} = O(e^{-K}). \eqno (3-30)$$
Let $Q$ be the component of $J\ssm U$
which does not contain $c$. Then 
$$P(S,Q)=O(\mu\lambda_*)=O(\lambda_*e^{-K}), $$
and hence
$$ {|\tilde U|\over |\tilde S|}= O(\lambda_*e^{-K}) \eqno (3-31) $$
as well. estimates (3-30) and (3-31) imply
$$ P(\tilde I, \tilde U)\geq 2 K+\log^+{1\over \lambda}-O(1). \eqno (3-32)$$
Pulling this back by the quadratic map, we get (3-23).\QED

\proclaim Lemma 3.10.  Let both intervals $I_s$ and $I_t$ be non-central 
and lie on the same side of $c$. Then

$$Q(G')\succ {1\over 2}K + \log{1\over \mu}-O(1) \eqno (3-33)$$
or
$$Q(G')\geq K+B(K,\lambda)-O(1) \eqno(3-34)$$
where $B(K,\lambda)\geq 0$, and $B(K,\lambda)=K/4-O(\lambda)$ for 
$\lambda\leq\bar\lambda$.

{\bf Proof.} Let us again consider several cases depending on the 
  non-linearity $\alpha$ of the quadratic map $\phi$ on the non-central 
 intervals. Let $\bar\alpha>0$ be small.

{\sl Case 1.} Let $p=1$ or $\alpha\leq\bar\alpha.$

\noindent For $p=1$ let us use representation (3-0). For $p>1$
the condition $\alpha\leq\bar\alpha$ holds.  Hence for any $p\geq 1$
$$g^p|T'=\psi\circ \phi|T' \eqno (3-36)$$ where $\psi$
is a diffeomorphism with bounded non-linearity.
Hence pulling $T$ back  by $\psi$, we don't spoil (3-16).
Composing this with $\phi$-pull-back, we get (3-33) by Lemma 3.8.  



{\sl Case 2.} Let $p>1$ and $\alpha\geq \bar\alpha$.
Then  Lemma 3.5 yields $\mu=O(e^{-K})$.

 Let $W\supset g T$ be a $g^{p-1}$-fold pull-back of $J$. Given an interval
$X\subset J$, denote by 
$\tilde X\subset W$ its $g^{p-1}$-fold pull-back. 
Pulling the pairs of intervals $I_s\supset V$ and
$I_t\supset U$ back to $W$ we get
$$P(\tilde V : \tilde I_s)=O(\mu)=O(e^{-K})\quad {\rm and}\quad
  P(\tilde U : \tilde I_t )=O(\mu)=O(e^{-K}). \eqno (3-37)$$
Let us pull the interval $I_s\cup H\cup I_t$ back subsequently by
$g$ and then by $g^{p-2}$. Apply Lemma 3.8 on the first step 
and the Koebe Principle on the second. This yields
$$P(\tilde I_s, \tilde I_t)\succ {1\over 2} P(H)\geq 2B(K,\lambda) \eqno (3-38)$$
where $B(K,\lambda)$ is as was claimed.
Estimates (3-37) and (3-38) yield
$$P(\tilde U, \tilde V)\geq 2 K+2 B(K,\lambda)-O(1).$$
Pulling this back by $g$, we obtain (3-34).
\QED 

\proclaim Lemma 3.11. Let $I_s$ and $I_t$ be non-central intervals 
lying on the opposite sides of $c$. Then
      $$Q(G')\geq K +B(K,\lambda)-O(1) $$
with $B(K,\lambda)$ as in Lemma 3.10.

\noindent {\bf Proof.} Let us again consider two cases.

{\sl Case 1.} Let $p=1$ or $\alpha<\bar\alpha$.

\noindent Then argue as in Case 1 of the previous lemma but use Lemma 3.7
  instead of Lemma 3.4. This yields $Q(G')\geq (5/4)K-O(1)$ which is better
  than what is claimed.

{\sl Case 2.} Let $p>1$ and $\alpha\geq \bar\alpha.$

\noindent  Then  argue as in Case 2 of the previous lemma.  \QED

\medskip\noindent {\bf More relations between the parameters.}
Let us mark the quantifiers  of level $n+1$ by ``prime": 
$\mu'\equiv\mu_{n+1},\; \lambda'\equiv\lambda_{n+1}$ etc.
 The following lemma provides us with rough estimates of the 
parameters of level $n+1$ through $\mu=\mu_n$.

\proclaim Lemma 3.12. The following estimates hold:
$$\lambda'=O(\sqrt{\mu}) \quad {\rm and} \quad \mu'=O(\sqrt{\mu}),
      \eqno (3-39)$$ 
$$\rho'=O(\mu),   \eqno (3-40)$$
$$K'\geq{1\over 2}\log{1\over \mu}-O(1). \eqno (3-41) $$
In the case of non-central return, (3-39) can be improved as follows:
$$ \mu'=O(\sqrt{\mu\lambda^*}). \eqno (3-42)$$

\noindent {\bf Proof.} Let us take an interval $U'$ of level $n+1$
and consider its image $U=g U'\subset I_t$.  If $t=0$ we have 
$P(U|J)\leq P(I_t|J)\approx 4\mu.$
Otherwise by Lemma 3.1  and estimate (3-2) we have
$P(U|J)=O( \lambda^*\mu).$
Now estimates (3-39) and (3-42) follow from decomposition (3-0).

It follows from (3-0) and Lemma 2.2 
that $g_{n+1} |U'=\psi\circ \phi|U'$ where
$\psi$ is a diffeomorphism with the Koebe space  spreading over
$J$. Since $\psi(\phi U')\subset I$, (3-40) follows.

In order to get (3-41) let us take two intervals $U'$ and $V'$ of level
$n+1$ and go through our basic construction (see the ``Strategy"). 
Represent $g^p|T'$ as a composition $\psi\circ \phi$
where $\psi$ is a diffeomorphism with a Koebe space spreading over
$J$. Now pull the pair of intervals $I_t\supset U$ back by $g^p$
taking into account that $P(U|I_t)\prec 2\mu$. We see that
$|U'|/|G'|=O(\sqrt{\mu})$, and the result follows.
                                                \QED 

In what follows we restore the index $n$. Let us now treat the problem of
estimating $\mu_{n+1}$ through the parameters of lower levels. 

\proclaim Lemma 3.13. Let $\epsilon>0$. In the non-central return case 
one of the following estimates holds:
$$\mu_{n+1}=O(\sqrt{\lambda_n^*}\exp(-(1-\epsilon)K_n/2)), \eqno (3-43)$$
or
$$\mu_{n+1}=O( \exp(-(1+\epsilon/2)K_n/2)), \eqno (3-44)$$
or
$$\mu_{n+1}=O\left(\sqrt{\mu_n\mu_{n-1}\lambda_{n-1}^*}\right). \eqno (3-45)$$

\noindent {\bf Proof.} Let us again consider several cases.

{\sl Case 1.} Let $\alpha_n>\exp(-\epsilon K_n).$
 Then by Lemma 3.5,
$$\mu_n=O(\exp(-(1-\epsilon)K_n)).$$ Using this and (3-42) we obtain (3-43).
\smallskip

{\sl Case 2.}  Let $\alpha_n\leq \exp(-\epsilon K_n).$
 Let $U=g_n I^{n+1}\subset I^n_t,\; t\neq 0$. Let $L$ be 
the component of $I^n_t\ssm U$ which is closer to $c$ and $R$ be 
the other component.\smallskip 

\noindent
 {\sl High return subcase.}  Then arguing as in Lemma 3.3 we see that
$$\eqalign{\log{1\over \mu_{n+1}^2}\succ 
\log\left (1+{{\rm dist}(U,c)\over |U|}\right )\geq
 \log\left (1+{{\rm dist}(I^n_t,c)\over |I^n_t|}\right )+
\log\left (1+{|L|\over |U|}\right )\geq
\cr\geq \left [\log{1\over 2\mu_n}+{1\over 2}\log\left (1+{{\rm dist}(I^n_t,c)
\over |I^n_t|}\right )\right ]
+{1\over 2}\log\left (1+{{\rm dist}(I^n_t,c)\over |I^n_t|}\right )-O(1)\geq\cr
\geq K_n+{1\over 2}\log{1\over\alpha_n}-O(1)\geq (1+\epsilon/2)K_n-O(1),\cr}$$
and (3-44) follows.\smallskip

\noindent
{\sl Low return subcase.} Let $b$ be the boundary point of $I^{n-1}$ lying on
the same side of $c$ as $I^n_t$. 
Set $X=[c,I^n_t],$ the convex hull of $c$ and $I^n_t$, $Y=[ I_t^n,b]$. 
 We need to refine the situation again.

\noindent (i) Let $|Y|^2\geq |X|\cdot |I^n_t|\exp(\epsilon K_n/2)$.
Then
$$\eqalign{\log{1\over \mu_{n+1}^2}\approx
\log\left (1+{{\rm dist}(U,b)\over |U|}\right )\geq
 \log\left ({|Y|\over |I^n_t|}\right )+
\log\left (1+{|R|\over |U|}\right )\geq\cr\geq
\left [ {1\over 2}\log\left ({|X|\over |I^n_t|}\right )+
\log{1\over 2\mu_n}\right ]+{1\over 2}\epsilon K_n -O(1)
\geq (1+\epsilon/2)K_n - O(1),\cr}$$
and we  have (3-44) again.\smallskip

\noindent (ii) Let $|Y|^2\leq |X|\cdot |I_t^n|\exp(\epsilon K_n/2)$.
Then ``an exponentially low return" occurs:

 $$|Y|/|I^{n-1}|\leq\exp(-\epsilon K_n/4).$$
It follows that
$${\rm dist}(g_{n-1}I^n_t,g_{n-1}b)/{\rm dist}(g_{n-1}I^n_t,g_{n-1}c)=
O(\exp(-\epsilon K_n/8)). \eqno (3-46)$$
Let $g_{n-1}I^n_t\subset I^{n-1}_j$. Then (3-46) implies that $I^{n-1}_j$
is a non-central interval, that is $j\neq 0$. 
Hence
$$|g_{n-1}I^n_t|/|g_{n-1}Y|=O(\mu_{n-1}\lambda_{n-1}^*).$$
Since $g_{n-1}|Y$ has distortion $O(\exp(-\epsilon K_n/4))$, 
we conclude that
$$|I^n_t|/|Y|=O(\mu_{n-1}\lambda_{n-1}^*)$$
as well. Together with $P(U|I^n_t)=O(\mu_n)$ this implies (3-45). \QED

Now we are prepared to prove the Conditional Version of Theorem B. 
To make life easier, let us first treat
the case when there are no central returns at all.\medskip

\noindent{\bf No central returns case.}
If $\mu_0$ is small then $K_1$ is big by (3-41). So, we can make the 
following inductive assumption: There is a $\theta>0$ such that
$$K_i\geq \theta i,\quad i=1,2,...,n, \eqno (A_n)$$
and 
$$\mu_i\leq {1\over 2} \exp(-\theta i/2),\quad i=1,2,...,n. \eqno (B_n)$$
By (3-39) we have
$$\lambda_i\leq\exp(-\theta i/4)<<\bar\lambda,\quad i=1,2,...,n. \eqno (3-47)$$
Now Lemma 3.13 allows us to conclude that there is a $\delta>0$ such that
$$\mu_i\leq\exp(-(\theta+\delta) i/2),\quad i=1,2,...,n+1. \eqno (3-48)$$
which is certainly stronger than $(B_{n+1})$.

In order to obtain $(A_{n+1})$ let us take a gap $G'$ between two intervals
$U'$ and $V'$ of level $n+1$ and push it forward as described in the above
``Strategy". Then we will find two intervals 
$I^n_s$ and $I^n_t$. Let us consider three cases depending on the position
of these intervals.\smallskip
{\sl Case 1. Let $I_s^n$ be the central interval.} Then by Lemma 3.9 and
estimate (3-47) we conclude that there is an absolute constant $\chi>0$
such that $Q(G')>K_n+\chi$ which is greater than $\theta(n+1)$, provided
$\theta$ was selected to be smaller than $\chi$. Taking the infimum over all
gaps $G'$ we obtain $A_{n+1}$.            \smallskip   

{\sl Case 2. Let $I^n_s$ and $I^n_t$ be two non-central intervals lying on the
opposite sides of $c$.} Then  Lemma 3.10, assumption $A_n$ and estimates (3-47),
(3-48) give us a small $\delta>0$ such that
$$Q(G')\geq(1+\delta)(n+1)\theta, \eqno (3-49)$$
which implies $A_{n+1}$.\smallskip

{\sl Case 3.  Let $I^n_s$ and $I^n_t$ be two non-central intervals lying on the
same side of $c$.} Then Lemma 3.11, the assumption $(A_n)$
and (3-47) yield (3-49) again. \bigskip

\noindent{\bf Cascades of central returns.}
Let us have a non-central return on  level $m-1$ followed by 
the cascade
of central returns on levels $m, m+1,...,m+q-1$, and completed  by a
non-central return on level $m+q,\quad q\geq 1$. 
So,  $m,\; m+q+1\in {\cal L}$.  Set $g=g_{m+1}|I^{m+1}$
(see Figure 5). Then
$$g(c)\in I^{m+q}\ssm I^{m+q+1} ,\;{\rm and}\;
 g_{m+i}|I^{m+i}=g|I^{m+i},\; i=1,...,q+1,  \eqno (3-50)$$
and $I^{m+i+1}$ is the $g$-pull-back of $I^{m+i},\; i=0,...,q$.
Let us call $q$  the {\sl length} of the cascade. 


\vskip 2in
\centerline{Figure 5}\medskip

Let us fix a big natural number $N$.  Let us define 
$\omega(m)$ in the following way. If a non-central return on
level $m-2$ occurs, that is $m-1\in{\cal L}$, then set $\omega(m)=0$.
Otherwise the level $m-2$ completes a cascade of central returns of length $p$.
Then set $\omega(m)=\min(p,N)$.

Let us assume by induction that there are  $\theta>0$ and $\delta>0$ such that
$$K_{i+1}\geq ((\kappa(i)+\omega(i))\theta,\qquad i\leq m,\; i\in {\cal L} 
\eqno(A_m),$$ and
$$\mu_{i+1}\leq\mu_1\exp(-(\theta+\delta)\kappa(i)/2), 
  \qquad i\leq m,\; i\in {\cal L}. \eqno(B_m)$$
Our goal is to check $(A_{m+q+1})$ and $(B_{m+q+1})$,
 provided $\theta$ and $\delta$ are small enough. 

When we travel along the cascade of central returns the trouble is
 that the scaling factors $\mu_{m+i}$ is definitely increasing 
(and very fast: as $(\mu_m)^{1/2^i}$).
However, they are still quite small  ($<\bar\mu$) in the initial segment of the
cascade, so that we can apply all above lemmas.
 If $\mu_1$ is small enough then $(B_m)$ guarantees that for $i\leq N$
$$\mu_{m+i}\leq \bar\mu. \eqno (3-51) $$
Moreover,  both $\lambda_n$ and $\alpha_n$ are exponentially small, that is
setting $\kappa=\kappa(m)$ we have
$$\lambda_{m+i}=O(\exp(-\kappa\theta/2))\quad {\rm and}\quad
\alpha_{m+i}=O(\exp(-\kappa\theta)/2)),\; i=2,...,q+1. \eqno (3-52)$$
Indeed, take a non-central interval $I^{m+i}_k$ and push it forward by
$g^{i-1}$. Since $g^{i-1}(c)\in I^{m+1}$ while 
$g^{i-1} I^{m+i}_k\subset I^m\ssm I^{m+1}$, there is a non-central
interval $I^{m+1}_t$ containing $g^{i-1} I^{m+i}_k$. Let $X\supset I^{m+i}_k$
 be the pull-back of $I^{m+1}_t$  by $g^{i-1}$. Then $X$ is contained in
 $I^{m+i-1}\ssm I^{m+i}$ and $P(I^{m+i}_k|X)=O(\mu_{m+1}).$
This estimate together with $(B_m)$ implies (3-52).

These considerations also show that $g_{m+i}$ can be represented as a 
composition
of the quadratic map $\phi$ and a diffeomorphism whose Koebe space is spread
over $I^m$. Hence, the distortion parameters remain small:
$$\rho_{m+i}=O(\exp(-\kappa\theta/2)),\quad i=2,...,q+1.\eqno (3-53)$$

\smallskip\noindent
{\sl An estimate for $K_{m+2}$}.
A  trouble with this estimate is  that $\lambda_{m+1}$ need not be small.
 However, by the
induction assumption and (3-39) the only way this can happen is if
$m-1\not\in{\cal L}$ and  $m-2$ completes a long cascade of central returns, 
that is $\omega(m)=N$ is big which makes the assumption $(A_m)$ stronger.  

More specifically, let us follow the above ``Strategy".
Take a gap $G'$ between two intervals of level $m+2$ and push it forward
by iterates of $g_{m+1}$ until its endpoints are separated  by different
intervals $I_s^{m+1}$ and $I_t^{m+1}$ of level $m+1$. 
As usual, let us consider
several cases depending on the positions of these intervals.

\smallskip\noindent
{\sl Case 1.} Let $I_s^{m+1}$ be central. Then as we have explained
either $\lambda_{m+1}<\bar\lambda$ or
$$K_{m+1}\geq \kappa\theta +A \eqno (3-54)$$
 with a big $A$.
In both cases  Lemma 3.9 yields 
$$Q(G')>(\kappa+1)\theta \eqno(3-55),$$
provided $\theta$ is small enough.

\smallskip 
\noindent
{\sl Case 2}. Let  $I_s^{m+1}$ and $I_t^{m+1}$ be non-central lying on the
   same side of $c$. If $\lambda_{m+1}<\bar\lambda$ then Assumptions $(A_m)$,
$(B_m)$ and Lemma 3.10 imply that there is an $\epsilon>0$ such that
$$ Q(G')>(1+\epsilon)\kappa\theta\eqno(3-56),$$
which is certainly better than (3-55).

If $\lambda\geq\bar\lambda$ then (3-54) holds. Together with Lemma 3.10
this yields (3-55).

\smallskip\noindent
{\sl Case 3.} Let $I_s^{m+1}$ and $I_t^{m+1}$ be non-central intervals 
  lying on the
   opposite sides of $c$. Then argue as in the previous case using
  Lemma 3.11 instead of 3.10.

\smallskip So, in all cases (3-55) holds, and hence 
$K_{m+2}\geq \theta (\kappa+1)$.

\medskip\noindent
{\sl Estimates for $K_{m+i+1},\; i>1$, in the initial segment of the cascade,
 (while (3-51) holds).}  Now $\lambda_{m+1}$ is exponentially 
small by (3-52) but $\mu_{m+i}$ need not be exponentially small.
Let us assume by induction that
$$K_{m+j+1}\geq (\kappa+j)\theta,\; j=1,...,i-1. \eqno (3-57)$$
To pass to the next level let us apply again our strategy and go through the 
same bunch of cases depending on the positions of $I_s^{m+i}$ and $I_t^{m+i}$.
Cases 1,2,3 mean the same as above. 

\smallskip\noindent {\sl Case 1.}
Then Lemma 3.9 and (3-52), (3-57) yield
$$Q(G')\geq (\kappa +i)\theta. \eqno (3-58)$$

\smallskip\noindent {\sl Case 2.}
 Let $H$ be the gap between $I_s^{m+i}$ and $I_t^{m+i}$.
If $|H|\geq{\rm dist (H,c)}$ then Lemmas 3.6 and 3.8 yield the desired
estimate.
Otherwise $g$ has a bounded distortion on $T'=U'\cup G'\cup V'$
(the notations are the same as in  the Strategy description),
and it follows from Lemma 3.4 that 
$$P(G')>K_{m+i}+2\log{1\over\bar\mu}-O(1)>K_{m+i}+\chi,$$
provided $\bar\mu$ is small enough.

\smallskip\noindent
{\sl Case 3} is treated  in the standard way using Lemma 3.11 and  (3-52).

Conclusion: (3-57) follows for $j=i$.

\smallskip\noindent
{\sl Distortion control in the tail of the cascade.} Let $M$
be the first moment for which $\mu_{m+M}\geq\bar\mu. $  Since
$\mu_{m+k+1}\approx \sqrt{\mu_{m+k}}$ for $k\leq q$,
$$ |I_{m+q+1}|\geq c |I_M|  \eqno (3-59)$$
(with $c\approx \bar\mu^2 $). 

\proclaim Lemma 3.14. The map $g^{ (j+1)}$ has  bounded
 distortion on both components of \break
$I^{M+j}\ssm I^{M+j+1}$, $0\leq j\leq q$. 

\noindent{\bf Proof.} As we know $g$ is almost quadratic. Hence by (3-59)
it has bounded non-linearity $n_g$ on $I^M\ssm I^{m+q+1}$. Let $L$ be a component
of $I^{M+j}\ssm I^{M+j+1}$. Then $f^k L$ is a component of
$I^{M+j-k}\ssm I^{M+j+1-k}$. Hence by  the standard argument
the non-linearity of $g^{j+1}$ on $L$ is bounded by
$$n_g \sum_{k=0}^j |g^k L|\leq n_g |I^M\ssm I^{m+q+1}|=O(1). $$  \QED 

Pulling now the intervals from level $M$ back to the tail of the cascade,
we conclude that
$$K_{M+i}\geq K_M-O(1)\geq (\kappa+M)\theta,\quad  i\leq q-M+1  \eqno (3-60)$$
(for, perhaps a bit smaller $\theta$).

\medskip\noindent {\bf A Markov scheme.} 
 Let us build up a Markov map $F$. Let $N$ be a big number as 
selected  above.
Set $K_j^{m+2N}=I_j^{m+2N}, \; j\neq 0$,
and pull these intervals back
by iterates of $g$ to levels $m+i,\; i=2N,...,m+q+1.$ Denote the corresponding
intervals by $K_j^{m+i}.$ Now set
$$F|K_j^{m+i}=g \quad {\rm for}\; 
i>2N\quad {\rm and}\quad F|K_j^{m+2N}=g_{m+2N}. $$
This map $F$ carries $K_j^{m+i}$ onto $K_j^{m+i-1}$ for $i>2N$,
 and carries $K_j^{m+2N}$ onto $I^{m+2N-1}$ covering all 
intervals of our scheme.

\medskip\noindent {\bf Proof of $(A_{m+q+1})$.}
Take two intervals $U'=I^{m+q+2}_k$ and $V'=I^{m+q+2}_l$ of level $m+q+2$,
consider their images by $g$, and 
then push them forward by iterates of 
$F$ until the first moment they don't belong to the same interval 
$K^{m+i}_j$ of our scheme. This long-term composition is almost quadratic
as one can see from (3-52). 
Denote the corresponding images  of $U'$ and $V'$
by $U$ and $V$.   We again have to consider several cases.

\smallskip {\sl Case 1.}  Assume that for some $N\leq j\leq q+1$ $V$
 belongs to a 
central interval   $I^{m+j}$ while $V\subset I^{m+j}_t,\; t\neq 0$.
If $j>M$ push these intervals forward to level $M$. By Lemma 3.14 this results
in a bounded change of the Poincar\'{e} length of the gap between the 
central interval $I^{m+j}$ and $V$. Let us denote new intervals by
$\tilde U$ and $\tilde V$. Now they lie on the  level $l=\min(M,j)$.
Then Lemma 3.3,  estimate (3-52) and the above estimates of $K_{m+i}$ yield
$$P(I^l,V)\geq 2 K_{m+2N}\geq 2( m+N)\theta.$$
 Pulling this back by the quadratic map postcomposed with a bounded distortion
 map, we get the desired estimate. 

\smallskip {\sl Case 2.} Let $V\subset I_s^{m+j}$, $U\subset I_t^{m+j}$ with
  $s\neq 0$ and $t\neq 0$. As in the previous case, pushing these intervals
  forward, we can assume that $j\geq M$. Let $H$ be the gap between 
  $I_s^{m+j}$ and $I_t^{m+j}$.

\smallskip {\sl Subcase (i).} Let $H\ni c$. The the standard argument based
  on Lemma 3.7 gives the desired estimate.

\smallskip {\sl Subcase (ii).} Let $H\not\ni c$ and $|H|\geq {\rm dist} (H,c)$.
   Use Lemma 3.6 instead of 3.7.

\smallskip {\sl Subcase (iii).} $H\not\ni c$ and $|H| < {\rm dist} (H,c)$. 
Then let us push the intervals forward by iterated $g$. Set $U_n= g^n U$ etc.
Assume there is a moment for which $|H_n|\geq {\rm dist} (H,c)$.
Then $g^n|U\cup H\cup V$ has a bounded distortion for the first such  moment.
Hence we can argue as in the previous subcase.

If there is no such a moment, then push the interval to the very beginning
of the cascade (to level $m$) and apply Lemma 3.4.

\medskip\noindent{\bf Proof of $(B_{m+q+1})$.}
We should estimate  $\mu_{m+q+2}$. Let us push $I^{m+q+2}$ forward to higher
levels: 
$$T^{l}=G_l I^{m+q+2}\equiv g_{l}\circ g_{l+1}\circ...\circ g_{m+q+1}I^{m+q+2}.$$
Let us stop on the highest level $l < M$ for which one of the following
properties hold:

\noindent (i) the map $G_l$ is not exponentially low in the sense of 
Lemma 3.13, Case 2-ii (that is $G_l I^{m+q+1}$ belongs to the
$\exp(-\epsilon K_{l-1})|I^{l-1}|$-neighborhood of $\partial I^{l-1}$), or

\noindent (ii) $l-1\in {\cal L}$ and $\kappa(l-1) = \kappa(m)-1$.
This means that we have arrived at the beginning of the previous cascade.

It follows from Lemma 3.14 that $G$ is a quadratic map postcomposed 
with a bounded distortion map. 
This allows us to apply Lemma 3.13 for $G_l$ instead of $g$, and
to estimate $\mu_{m+q+2}$ through the parameters of level $l$.
If (i) occurs then (3-43) or (3-44)-like estimates hold.
Together with the above estimates of $K_l$ and $\lambda_l$ they
yield the desired estimate.

Otherwise there is a non-central interval $I^l_s\supset T_l$ which can
be monotonically pulled back by $G_l$. Since
$$P(T_l|I^l_s)=O(\mu_l\mu_{l-1}...\mu_{m+q+1})=O(\exp(-5\theta/4)),$$
we conclude that
$$\mu_{m+q+2}=O\left(\sqrt{P(T_l|I^l_s)}\right)=O(\exp(-5\theta/8)),$$
and $(B_{m+q+1})$ follows.

\bigskip
\centerline{\bf \S 4. Real bounds and limits of renormalized maps.}\medskip
In this section we will prove Theorem B for maps with a non-minimal critical
set $\omega(c)$ and for maps of unbounded type. For maps of bounded type
 will show that if the scaling factors stay away
from zero then the family of renormalized maps is compact, and all limit
maps are real analytic. As usual, let us
first assume that {\sl there are no central returns}.
Set $x_m=f^m x$.
\medskip

\noindent {\bf  A priori bounds.} As in \S 3 let us consider
the decomposition
(3-0) $g_n|I^n_i=h_{n,i}\circ\phi$, and denote by $\rho_n$ 
the maximal distortion of the diffeomorphisms $h_{n,i}$. 

\proclaim Theorem 4.1 (Martens [Ma]). The distortions $\rho_n$ are
uniformly bounded.

\proclaim Lemma 4.2. The scaling factors $\mu_n=|I^n|/|I^{n-1}|$
are  bounded away from 1.

\noindent
{\bf Proof.} The case of low return on level $n-1$ was treated in
   ([Ma], Lemma 3.7).

In the case of  high return let us assume that $\mu_n$ is close to 1.
Then because of the bounds given by Theorem 4.1, the next scaling factor
$\mu_{n+1}=O(\sqrt{1-\mu_n})$ will be very small. 
Then  by the results of \S3 $\mu_m\to 0$. \QED  

\proclaim Lemma 4.3. All Poincar\'{e} lengths $P(I^n_j|I^{n-1})$
  are bounded away from $\infty$.

\noindent {\bf Proof.}  This follows from the previous lemma and Lemma 3.12.
                                                                  \QED

Let $\bar\lambda$ denote  an upper bound of Poincar\'{e} lengths
$P (I^n_j | I^{n-1})$.
Consider two intervals $T\supset G$ with $L$ and $R$ to be the components of
$T\ssm G$. Denote by
 $$\bar\sigma={e^{\bar\lambda}\over 1+e^{\bar\lambda}}<1 \eqno (4-1)$$ 
an upper bound of  
$|G|/(|G|+|L|)$ provided $P(G|T)\leq\bar\lambda$.

\smallskip\noindent
{\bf Orders and ranks.} Let us define the {\sl order} 
ord$^n_0$ of the return
to level $n$ as the return time of $g_n$-orb($c$) back to $I^n$. Let us
also define $l$-orders ord$^n_l$ as the return time of $g_{n-l}$-orb($c$) 
back to $I^n$. In terms of the return graph the ord$^n_0$ is just the number
of edges beginning at $I^{n+1}$ (and leading to the previous level $n$).
The ord$^n_l$ is the number of paths of length $l+1$ beginning at $I^{n+1}$
(and leading to the level $n-l$).

\proclaim Lemma 4.4. If the scaling factors $\mu_n$ stay away from 0 
then for each $l$ the $l$-orders ord$^n_l$ of returns to all levels 
are uniformly bounded.

\noindent {\bf Proof.} If on level $n$ a return of high order $p$ occurs then
by Lemma 4.3 the next scaling factor $\mu_{n+1}=O(\sqrt{\bar\lambda^p})$ 
is very small.
Similarly, for a given $l$, if ord$^n_l$ is big then traveling down the
graph from level $n-l$ to $n+1$ we see that $\mu_{n+1}$ is small as well.
\QED 

Let us consider now the Markov family ${\cal M}$  of intervals 
obtained by pull-backs
of the initial interval $I^0$ (see \S 2).

Let us assign to the  critical intervals $K\in {\cal M}$  rank 0.
Let us say that an interval  $K\in {\cal M},\; K\subset I^{n-1}\ssm I^n$
has {\sl rank} $k,\; k\geq 1$, if orb($c)$ passes through
it before the first return to $I^{n+k-1}$ but after the first return to
$I^{n+k-2}$. For example, $k=1$ if orb($c$) passes through $K$ before
the first return to $I^n$.
For $K=I^n_j$ this can be nicely expressed in terms of the return graph as
the length $k$ of the shortest path
leading from $I^n_j$ down to a central interval $I^{n+k}$.

\proclaim Lemma 4.5. Let $K\in {\cal M}, \; K\subset I^{n-1}\ssm I^n$,
 and rank$(K)\geq k>0$. Let us take a point $x\in I^{n+k-2}\ssm I^{n+k-1}$,
 and consider the first moment $l$ when $f^l x\in K$ (provided there is one).
 Then the interval $K$ can be
diffeomorphically pulled back along the orb$_l(x)$ to an interval
$K'\subset I^{n+k-2}\ssm I^{n+k-1}$.

\noindent {\bf Remark.} We don't claim that $K'\ni x$ but just some iterate
$x_s,\; s\leq l,$ so that the pull-back has length $l-s$.\smallskip 

\noindent
{\bf Proof.} 
 If $k=1$  there is nothing to prove (set $K'=K$). If $k>1$ let us
consider the last moment $s < l$ when the orb$(x)$ visits $I^{n}$.
Then there is a moment $p$ such that $g_n^{\circ p}(x_s)\in K$ and all
intermediate iterates $g_n^{\circ m}(x_s)\not\in I^n,\; 0<m<p$. I claim that
the pull-back of $K$ along the $g_n$-orbit of $x_s$ is monotone.
Indeed, otherwise $g_n^{\circ p}(c)\in K$ while 
$g_n^{\circ m}(c)\not\in I^n,\; 0<m<p$,
so that rank$(K)=1$.

Let $K_1\ni x_s$ be the monotone pull-back of $K$ along the $g_n$-orbit of
$x_s$. Then $K_1\subset I^n\ssm I^{n+1}$ and rank$(K_1)\geq k-1$.
So, we can proceed by induction.\QED 

\proclaim Lemma 4.6. If the scaling factors $\mu_n$ stay away from 0 then
ranks of all intervals $I^n_j$ are uniformly bounded.

\noindent
{\bf Proof.} Let rank$I^n_j=k$.
Let $l$ be the first moment when orb($c$) visits $I^n_j$.
Then by the definition of rank there is an $s<l$ such that 
$c_s\in I^{n+k-2}\ssm I^{n+k-1}$.
By the previous lemma we can monotonically pull
$I^n_j$ to the level $n+k-2$  along the  orb$c_s$. We obtain an interval
$K'\subset I^{n+k-2}\ssm I^{n+k-1} $.

Let us consider the interval $I^{n+k-1}_i$ containing $c_s$.
Then one can see by induction (involving the Schwarz lemma) that 
$$P(I^{n+k-1}_i | I^{n+k-2})\leq P(I^{n+k-1}_i | K') \leq \bar\sigma^{k-1}$$
with $\bar\sigma$ from (4-1).
 Since orb($c$) passes through $I^{n+k-1}_i$ before its return to 
$I^{n+k-1}$, the scaling factor $\mu_{n+k}=O(\bar\sigma^{(k-1)/2})$
 is small if $k$ is big. 
 Contradiction.       \QED

\medskip\noindent {\bf Non-minimal case.} Now we are ready to proof Theorem B
in the non-minimal case.

 \proclaim Lemma 4.7. Assume that the critical set $\omega(c)$ is not minimal.
  Then the scaling factors $\mu_n$ go down to 0, $n\in {\cal L}$.
   
\noindent{\bf Proof.} Indeed, in the non-minimal case
 there is a level $n$ and a point $x\in \omega(c)\cap I^{n-1}$
which never passes through the central interval $I^n$.  It follows that the
return time of points of orb$(c)$ back to $I^n$ is unbounded and hence
there are infinitely many intervals $I^n_k$ of level $n$. The ranks of these
intervals certainly must grow up to $\infty$. Now lemma 4.6 provides us with
the starting condition for Theorem 3.0. \QED
   
\medskip\noindent
{\bf Unbounded Combinatorics.}  Assume that $\omega(c)$ is minimal.
Let us say that $f$ is {\sl of bounded type} if the number of intervals on
all levels is uniformly bounded, and of {\sl unbounded type} otherwise.
The unbounded case can also be treated by a purely real argument
(I owe this remark to Swiatek).

\proclaim Lemma 4.8. If $f$ has unbounded combinatorics then
the scaling factors $\mu_n$ go down to 0, $n\in {\cal L}$.

\noindent {\bf Proof.} If the combinatorics are unbounded then either
the $l$-orders of returns or ranks of the intervals are unbounded
(consider the return graph from \S 2). Now the required starting condition
for Theorem 3.0 follows from lemmas 4.4. and 4.6.  \QED  

\medskip\noindent
{\bf Bounded Combinatorics.} This is the main case when we need to involve
complex analytic methods. In  this subsection we will show that 
provided the scaling factors stay away from 0, there is
a sequence of renormalized maps $C^1$-converging
to an analytic  map. 

Here by the gaps of level $n$ we will mean the components of
$I^{n-1}\ssm\cup I^n_j$.

\proclaim Lemma 4.9. If the scaling factors stay away from 0  then 
all intervals and all gaps of level $n$ are commensurable
with $I^{n-1}$.

\noindent {\bf Proof.} Let us show first that the intervals $I^n_i$ may not
be tiny as compared with $I^{n-1}$. Indeed, because of Lemma 4.3 such an 
interval should lie very close to $\partial I^{n-1}$. On the other hand,
it follows from Theorem 4.1 that $g_{n-1}$ is quasi-symmetric. Hence,
$g_{n-1} I^n_j$ lies very close to $\partial I^{n-2}$ and, moreover,
$g_{n-1} I^{n-1}$ covers the interval $I^{n-1}_s\supset g_{n-1} I^n_j$
(again because of Lemma 4.3). Hence we can monotonically pull $I^{n-1}_s$
back by $g_{n-1}$. 

Now we can apply the same argument to the interval $I^{n-1}_s$ and map 
$g_{n-2}$ and so on. In such a way we will find a big $l$  and an interval
$I^{n-l}_t$ which can be monotonically pulled back by
$g_{n-l}\circ...\circ   g_{n-1}$ along the orbit of $I^n_i$. Let
 $ K\subset I^{n-1}$ be this pull-back. This  interval provides us with
a big space around $I^n_i$, namely   $|I^n_i|/|K|\leq \bar\sigma ^l$.
Since rank$I^n_i$ is bounded by Lemma 4.6, we can pull this space back and 
obtain a small scaling factor. Contradiction.

Let us now consider the gaps. They may not be too big as compared with the
intervals since otherwise the intervals would be tiny as compared with
$I^{n-1}$. Let us consider any gap $G$ in between $I^n_i$ and $I^n_j$.
Arguing as in the proof of estimate (3-41) one can see that the Poincar\'{e}
length $P(G)$ is bounded from below provided the scaling factors are 
bounded away from 1.  This means that $G$ is not tiny as compared with
one of the intervals $I^n_j$, $I^n_i$. Since these intervals are commensurable
with $I^{n-1}$, so the interval $G$ is as well.

As to the two  ``boundary" gaps, they are not too small as compared with
the attached intervals because of Lemma 4.3.  \QED

Consequently, we can select a sequence of renormalized maps $R^n f$ 
in such a way
that the configurations of intervals and gaps converge to a non-degenerate
configuration of intervals and gaps. Use now the rescaled representation 
(3-0) for these maps: 
$$R^n f\; |\;\tilde I^n_i=G_{n,i}\circ\phi $$
where $\tilde I^n_i$ are the rescaled intervals of level $n$ and $G_{n,i}$
are diffeomorphisms of appropriate intervals onto the unit interval.
 By Theorem 4.1
and the previous lemma, the inverse maps $G_{n,i}^{-1}$ have uniformly
bounded $C^2$-norms, and hence form a $C^1$-compact family.  
So, we can select a $C^1$-convergent sequence of renormalized maps,
$$G_{n(s),i}^{-1}\to G_i^{-1}\quad {\rm as} \quad n(s)\to\infty .$$

Each $G_i$ is a long composition of the square root maps and 
diffeomorphisms whose total distortion is controlled by 
$$\omega_n=\sum_j\sum_{m=1 }^{p(n,j)-1}|f^m (I^n_j)|$$
where $p(n,j)$ is the return time of $I^n_j$ back to $I^{n-1}$. 
But $\omega_n\to |\omega(c)|=0$ by [Ma]. Hence the total distortion
of the diffeomorphisms involved is vanishing.

Now the  ``Shuffling Lemma"  
(see [S] or [MS], Ch. VI,Theorem 2.3) yields that the limit
of renormalized maps is real analytic and, moreover, belongs to so called
Epstein class which we are going to study in the next section.

\medskip\noindent {\bf Cascades of central returns.}
Let us have (as in the end of \S 3)  a non-central return on level $m-1$ 
followed by the cascade of central returns on levels
 $m, m=1,..., m+q-1$, and completed
by a non-central return on level $m+q$. The following remarks allow to
adjust the previous analysis to this case.

First of all, the first scaling factor $\mu_m$  of the cascade stays away
from 1 by lemma 4.2.

Assume now that the starting conditions don't hold, that is, 
 the scaling factors stay away from 0. Then $\mu_{m+q}$ stay away from 1.
For, otherwise we have a high return on the level $m+q-1$, and the next
scaling factor $\mu_{m+q+1}$ is tiny (see the argument of Lemma 4.2).

Furthermore, the ratio
$${|I^{m+q}|\over |I^m|}\geq\delta>0 \eqno  (4-2)$$
also stays away from 0. For, otherwise the non-central intervals $I_k^{m+q}$
have a small Poincar\'{e} length in the appropriate component of
$I^{m+q-1}\ssm I^{m+q}$. This would enforce $\mu_{m+q+1}$ to be small again. 

Consequently the map
$$g_m^{\circ (q-1)}: I^{m+q-1}\ssm I^{m+q}\rightarrow I^m\ssm I^{m-1}$$
has a bounded distortion. Indeed, since $g_m$ is  quadratic  up to
a bounded distortion, by (4-2) it has bounded non-linearity on 
$I^m\ssm I^{m+q}$. Since the iterates of $I^{m+q-1}\ssm I^{m+q}$ are disjoint,
the claim follows.

Now we can consider the return graph skipping all intermediate levels
between $m+q-1$ and $m$, and to define the orders and ranks of the intervals
through this graph. As we have shown, on all levels of the graph we have 
a priori bounds of the scaling factors, and passage from one level to
another has a bounded distortion. Now we can repeat the above argument.


\bigskip
\centerline{\bf \S 5. Epstein class and complex bounds.}\smallskip 
The goal of this section is to show that an appropriate renormalization
of an  analytic map  of Epstein class is polynomial-like. 

\medskip\noindent {\bf Polynomial-like maps.}
By a {\sl polynomial like map} we mean an analytic branched covering 
$$f: \cup_{i=0}^l U_i\rightarrow V$$ where $U_i$ and $V$ are topological disks,
cl$U_i\subset V$.  
Let us consider a class $\cal A$ of polynomial-like maps
$f$ having a single non-degenerate critical point $c\in U_0$.

\medskip\noindent {\bf Epstein class.}
Given an interval $I\subset {\bf R}$, let 
$P(I)={\rm int}\;({\bf C}\ssm ({{\bf R}\ssm I}))$
denote the complex plane slitted along two rays, $D(I)$ denote the disk
based upon $I$ as a diameter. Let $g$ be a real analytic map
$\cup I_k\rightarrow J$ satisfying the following properties:

{\item (i).} For $k\neq 0$ there is the inverse map $(g|I_k)^{-1}$ which
univalently maps $P(J)$ onto $P(I_k)$. 

{\item (ii).} $g|I_0=h\circ\phi $, $\phi(z)=(z-c)^2$ 
is the quadratic map, and
$h^{-1}$ has a univalent analytic continuation to $P(J)$.

Let us call this class of maps {\sl Epstein class} ${\cal E}$ 
(compare [S]).
Let us start with a general lemma from hyperbolic geometry
which was an ingredient of  Sullivan's Sector Lemma.
It is essential for our complex bounds as well.

\proclaim Lemma 5.1. 
Let $\phi: P(I)\rightarrow P(J)$ be an analytic map which
 maps $I$ diffeomorphically onto $J$. Then $\phi D(I)\subset D(J).$

\noindent
{\bf Proof.} The interval $I$ is a Poincar\'{e} geodesic in $P(I)$, and
the disk $D(I)$ is its Poincar\'{e} neighborhood (of radius independent
of $I$). Since $\phi$  contracts the  Poincar\'{e} metric,
 we are done. \QED

\noindent {\bf High returns.}
In order to make the following discussion more comprehensible let us 
dwell first on the case when
high returns occur on all levels (compare [LM], \S 8). 

\proclaim Lemma 5.2. Let $f\in {\cal E}$. Assume that we have high returns
on all levels.
Then $g_n$ is polynomial-like for some $n$. 

\noindent
{\bf Proof.} We can assume that the scaling factors stay away from zero.
Then given an arbitrary small $\sigma>0$, we can select a moment $n$
such that
$$\mu_n\geq \mu_{n-1}(1-\sigma). \eqno (5-1) $$
Set $q=g_n c\in I^{n-1}\ssm I^n,\; I^m=[\alpha_m,\beta_m]$ where $\beta_m$
lies on the same side of $c$ as $q$ (see Figure 6).


 Let us estimate
the Poincar\'{e} length $P_0$ of 
$[q,\beta_{n-1}]$ in $[\alpha_{n-1},\beta_{n-2}]$:
$$\eqalign{P_0\leq \log\left(1+{1-\mu_n\over 1+\mu_n}\right)+
\log\left(1+{\mu_{n-1}(1-\mu_n)\over 1-\mu_{n-1}}\right)
\leq \cr\log{2\over 1+\mu_n}
+\log(1+\mu_n)+o(1)=\log 2+o(1)\quad {\rm as}\quad \sigma\to 0\cr}
 \eqno (5-2)$$
(we have replaced $\mu_{n-1}/(1-\mu_{n-1})$ by  $\mu_n/(1-\mu_n)$ using
(5-1) and monotonicity of the linear-fractional function 
$x/(1-x),\; x\in (0,1)$,
and boundedness of $\mu_n$ away form 1 (Lemma 4.2)).
Let us consider now the representation  $g_n=h\circ\phi$ where
$\phi$ is the quadratic map and $h$ is a diffeomorphism of an interval
$T\supset \phi I^n$ onto $[\alpha_{n-1}, \beta_{n-2}]$. 
Set $T=\phi I^n\cup M\cup R$ where $M$ and $R$ are $h$-pull-backs of
$[q,\beta_{n-1}]$ and $[\beta_{n-1},\beta_{n-2}]$ respectively.
By the Schwarz lemma and (5-2)
$$\log\left( 1+{|M|\over |\phi I^n|}\right)\leq P_0\leq\log 2+o(1),$$
hence
$${|M|\over |\phi I^n|}\leq 1+o(1) \quad {\rm as}\quad \sigma\to 0. 
\eqno (5-3)$$

\vskip 2.5in
\centerline{Figure 6}\medskip

Let us now take the set $V=D(I^{n-1})$ and pull it back by $h$. By Lemma 5.1
we will obtain a domain $U'$ contained in $D(\phi I^n\cup M)$. Pulling this
domain back by the quadratic map $\phi$ we obtain by (5-3) a 
convex domain $U_0$ which is almost contained in the disk $D(I^n)$. 

On the other hand, pulling $V$ back
by the univalent branches of 
$g_n$, we will get by Lemma 5.1 domains $ U_k$ such that 
$U_k\subset D(I_k^n),\quad k\neq 0$.
We conclude  that the sets cl($U_k$)  are pairwise
disjoint and are contained in $V$. \QED

\medskip\noindent {\bf Cut-off iterates.}
Let us define now a ``cut-off" iterate $g_n^{\circ l} K$ of an interval $K\ni c$
inductively in the following way: 
$$g_n^{\circ l} K = g_n(g_n^{\circ(l-1)} K\cap I^j_n)$$
where $I^n_j\ni g_n^{\circ(l-1)} c$.  If $K=I^n$ then one boundary point of
its cut-off iterates belongs to
the $\partial I^{n-1}$. Let us select the first moment $l$ for which
$g_n^{\circ l}( I^n,c)\cap I^n\neq\emptyset$. In the low return case $l>1$.

\medskip\noindent {\bf Low returns (a particular case).} 
 
\proclaim Lemma 5.3. Let us select a level $n$ for which (5-1) holds, and find
an $l$ as described above. Assume that $g_n^{\circ l} c\not\in I^n$. 
Then  appropriate
pull-backs of $D(I^{n-1})$ form a\break polynomial-like map. 

\noindent {\bf Proof.} Denote $q=g_n^{\circ l} c$, $p=g_n^{\circ(l-1)} c$,
 $I^m= [\alpha_m, \beta_m]$ where $\beta_m$ lies on the same side of $c$ as
$q$. Let $I^n_i$ be the interval of level $n$ containing $p$,
 and let $p$ divides it into the intervals
 $L$ and $G$ with $G$ closer to $c$ than $L$. Note that
$g_n L=[\alpha_{n-1}, q]$, $g_n G=[q,\beta_{n-1}]$.

Let us consider representation (3-0). Since $h_{n,i}$ has a Koebe space spread
over $I^{n-2}$, one can get the following estimate in the same way as (5-3):
$${|\phi G|\over |\phi L|} \leq 1+o(1)\quad {\rm as}\quad \sigma\to 0.
                   \eqno (5-4)$$ 
Set $\mu_n=\mu$.
Since the non-linearity of $\phi$ on $I^n_i$ is at most $\log(1/\mu)$,
we obtain the following estimate:

$$ {|G|\over |L|}\leq {1\over\mu}\;(1+o(1)) \quad {\rm as}\quad \sigma\to 0.
 \eqno (5-5)$$

Let $R$ be the component of $I^{n-1}\ssm I^n_i$ containing the critical point
$c$. Let us estimate now the Poincar\'{e} length $P(L,R)$. It follows from (5-5)
that
$$|G|\leq {|I^n_i|\over 1+\mu}\; (1+o(1))
 \leq {1-\mu\over 1+\mu}\; |I^{n-1}|\; (1+o(1)).
               \eqno (5-6)$$
Hence
$$\eqalign{P(L,R)\leq \log\left(1+{|G|\over |L|}\right)+
\log\left(1+{|G|\over |R|}\right)+o(1)\leq\cr\leq
\log\left(1+{1\over\mu}\right)+\log\left(1+{1-\mu\over (1+\mu)^2}\right)+o(1)=
 \cr =\log\left(1+{2\over \mu^2+\mu}\right)+o(1) <
 \log\left(1+{1\over\mu^2}\right)+o(1).\cr}
                 \eqno (5-7)$$

Let $T\ni c$ be the pull-back of $I^{n-1}$ by $g_n^l$. Then 
$g_n^{l-1} T= L$. Decompose $g_n^{l-1}|T=h\circ\phi$ where 
$$h: \phi T\cup G'\cup R'\rightarrow L\cup G\cup R=I^{n-1}$$ 
is a diffeomorphism.
By the Schwarz lemma and (5-7),
$${|G'|\over |\phi T|}\leq {1\over \mu^2}\; (1+o(1)).  $$

Set $\nu=|T|/|I^{n-1}|.$ It follows from the a priori bounds of \S 4 
that $\nu/\mu=|T|/|I^n|$ stays away from 1. Hence the last estimate can be
rewritten as 
$${|G'|\over |\phi T|}\leq {\tau^2\over \nu^2} \eqno (5-8) $$
with an absolute constant $\tau<1$.

Let us consider now the disk $D(I^{n-1})$. By Lemma 5.1 its pull-back $V_0'$
 by $h$ is contained in the disk $D(\phi T)$. Pulling $V'$
 by the quadratic map $\phi$ we obtain a domain $V_0$ based upon the
 interval $T$.
 By (5-8) the outer radius of $V_0$ around $c$ is less than 
$(\tau/\nu)|T|=\tau |I^{n-1}|$. Hence cl$(V_0)\subset D(I^{n-1})$.

We have completed the construction of the central domain $V_0$. Let
 us now construct non-central domains of definition of our polynomial-like map.
First of all, take a non-central interval $I^n_i$ of level $n$ and pull the
disk $D(I^{n-1})$ back by the corresponding branch of $g_n^{-1}$.  By Lemma
5.1 we obtain a domain $W^n_i\subset I^n_i$. 

Let now $x\in I^n\ssm T$. Then there is a moment $s<l$ when $g_n$-orb$(x)$
is separated from $g_n$-orbit of $T$ by intervals of level $n$. In other
words $g_n^{\circ s} x\in I^n_j$ while $g_n^{\circ s} T\cap I^n_j=\emptyset$.
 Hence we can monotonically pull $I^n_j$ back to $x$ by
 $g_n^{\circ s}$ and obtain an
interval $T(x)$. Moreover, we can univalently pull the domain $W^n_j$
back to $x$ and obtain a domain $W(x)\subset D(T(x))$  based upon $T(x)$.
A map $g_n^{\circ (s+1)}$ univalently carries this domain onto $D^{n-1}$.

Let us consider the whole (finite) bunch of intervals $I^n_i$ and $T(x)$,
and the corresponding bunch of domains $W^n_i$ and $W(x)$.  Let
us redenote them as  $T_i$ and $V_i,\; i=1,...$ respectively. 
Since $V_i\subset D(T_i)$, these domains are pairwise disjoint. They are
also disjoint from the central domain $V_0$. Indeed, the domains
$\phi V_i,\; i=0,1,...$ are pairwise disjoint since they are contained in
the disks  based upon disjoint intervals.

So, we have a polynomial-like map $H:\cup V_i\rightarrow D(I^{n-1})$.
\QED 

\medskip\noindent
{\bf A remark on scaling factors.}
The map $H$ constructed above satisfies the following property:
the $g_n$-first return map
to $T_0$ coincides with the $H$-first return map to $T_0$. Indeed, if
$H|T_i=g_n^{\circ m},\; i\neq 0,$ then by looking through the construction we 
see that $g_n^{\circ k} T_i$ does not intersect $T_0$. 

Let $T_0\equiv T^0\supset T^1\supset...$ be the sequence of the central
intervals of the renormalized maps $R^{\circ N} H$, and $\nu_n$ be the
corresponding sequence of the scaling factors. Then the above property
of the first return maps yield that there is an $N$ such that
           $$  I^{n+N}\subset  T^n \subset I^{n+N-1}.$$ 
Hence the scaling factors $\nu_n$ can be estimated through the $\mu_n$
and vice versa. It follow that $\nu_n\to 0$ if and only if $\mu_n\to 0$.

\medskip\noindent
{\bf Low returns: a general construction.}
 The particular construction described above will be one step of the general 
construction. Let us start with a map 
$$h_1 : \cup T^{1,0}_i\rightarrow T^0$$ of class ${\cal A}$.
Let us order pairs $(n,j)$ of integer numbers lexicographically.
We will construct a  finite hierarchical family of intervals 
$T^{n,j}_i,\; (n,j)\in E\subset {\bf Z}^2,$ 
satisfying the following properties:

\smallskip\noindent
0) $E=\{(n,j): \; 0\leq n\leq N,\; 0\leq j\leq j(n)\}\;$, $T^0\equiv T^{0,0}_0$, 
   $T^{n,j}_0\equiv T^{n,j}$ are symmetric intervals containing $c$;

Let $(n',j')\in E$ be the lexicographic successor of $(n,j)\in E$. Then

\smallskip\noindent
1) $T^{n,j}\supset {\rm cl}\;( T^{n',j'}_i)$;

\smallskip\noindent
2)  There is a map  $h_n: \cup T^{n,0}_i\rightarrow T^0$ of  class
   ${\cal A}$ induced by $h_1$;

\smallskip\noindent
3) For $j>0$ the intervals $T^{n,j}_i$ are obtained from $T^{n,0}_i$
   by pulling back by the central branch of $h_n$. Moreover,
   $h_n c\in T^{n,j}$ for all $j$ except for the last one $j(n)$
   (the reader can recognize here the cascades of central returns).

\smallskip\noindent
4) For $n<N$, $h_n T^{n,0}\not\ni c $, while $h_N T^{N,0}\ni c$.

Under such circumstances we will also consider the map
$$H_n:\left(\bigcup_{m< n,\; i\neq 0}\; T^{m,j}_i\right)\cup T^n_0
    \rightarrow T^0$$ 
such that
$$H_n|\; T^{m,j}_i=h_m.$$
Note that $H_n$ is a map of class ${\cal A}$ with the non-escaping critical
point.

Assuming we have already constructed all intervals and maps ut to the level
$T^{n,0}$, let us do the next step. 

\smallskip
{\sl High case}. If $h_n T^{n,0}\supset T^{n,0}$  we stop.

\smallskip
{\sl Central-high case.} If $h_n c\in T^{n,0}$ and $h_n T^{n,0}\ni c$,
we consider the cascade of central returns until the first moment of high
return. It produces the intervals $T^{n,j}_i$ by pulling 
$T^{n,0}_i$ back by the central branch of $h_n$. Then we stop.

\smallskip {\sl Low case.} Acting as in the above particular case
let us consider cut-off iterates $H_n^{\circ m} T^{n,0}$ 
of the central interval 
until the first moment $l$ when 
$$H_n^{\circ l}\;  T^{n,0}\cap T^{n,0}\not=\emptyset. $$
Then let us pull $T^0$ back by $H_n^{\circ l}$.
It gives us the central interval $T^{n+1,0}$. 

Similarly we will construct a non-central interval 
$T^{n+1,j}(x)\equiv T^{n+1,j}_i, \; x\in \omega(c)$, as 
the pull-back of appropriate $T^{m,j}_i,\; m< n$, corresponding to
the first moment $k$ when the $orb(c)$ is separated from the $orb(x)$:
$$H_n^{\circ k} x\in T^{m,j}_i,\quad H_n^{\circ k} c\not\in T^{m,j}_i.$$
Define $h_{n+1}|T^{n+1,j}(x)=H_n^{\circ (k+1)}$.

\smallskip {\sl Central-low case.}
Let $h_n  T^{n,0} c\in T^{n,0}$ and $h_n T^{n,0}\not\ni c$. 
Then let us consider the cascade of central
returns until the first one which is low. 
 It produces the intervals $T^{n,j}_i$ by pulling
$T^{n,0}_i$  and $T^{m,j}_i,\;( m<n, i\neq 0)$
 back by the central branch of $h_n$.
Now let us define a  map 

$$F: \left(\bigcup_{m\leq n,\; i\neq 0} T^{m,j}_i\right)\; \cup  T^{n,j(n)}_0
     \rightarrow T^0. $$ 

\noindent For $m<n$ set $F|\; T^{m,j}_i=H_n $.

\noindent For $m=n, i\neq 0$, set $F|T^{n,j}_i=H_n^{\circ (j+1)} $.

\noindent Finally set $F|T^{n,j(n)}_0=H_n$.

Now taking cut-off iterates $F^{\circ m} T^{n,j(n)}$ we can construct
the intervals $T^{n+1,j}_i$ in the same way as in the low case. 

\noindent
\medskip {\bf Class} $\tilde{\cal T}$. Observe that the above map $F$
does not belong to class ${\cal T}$: the image of the central interval
belongs to int$T^0$. Such a situation always occurs in the end of a cascade 
of central returns. In order to handle it we need to introduce a wider
class of maps.

Let us consider an interval $T^{-1}$ whose interior contains finitely
many disjoint closed intervals $T^0_i$, with $T^0\equiv T^0_0$ containing 
 $c$ and symmetric with respect to it.
Let $g: \cup T^0_i\rightarrow T^{-1}$ be a map with negative Schwarzian
derivative  and a single critical point $c$ satisfying the following
properties:

\smallskip\noindent
1) $g$ diffeomorphically maps any non-central interval $T^0_i,\; i\neq 0$
onto $T^{-1}$; 

\smallskip\noindent
2) $g|T^0$ is a composition of a quadratic map and a diffeomorphism
   onto $T^{-1}$  with negative Schwarzian derivative; 

\smallskip\noindent
3) If $g(T^0)\cap T^0_i\not=\emptyset$ and $g c\not\in T^1_i$ then
     $g(T^0)\supset T^1_i$ (``Markov property").

\smallskip\noindent 
4) $g c\not\in T^0$ (non-central return).

\smallskip\noindent  Let $\tilde{\cal T}$ be the class of such maps.

Observe that the first renormalization  of a map 
of class $\tilde {\cal T}$ belongs to class ${\cal T}$.
Moreover, if the scaling factor $\mu=|T^0|/|T^{-1}|$
is small then the scaling factor of the renormalized map is also small.
It follows that the results of \S 3 are still valid if we start with a
map of class $\tilde{\cal T}$: if $\mu<\delta$ then the scaling factors
of renormalized maps go down to 0. Let us find the best $\delta$
satisfying this property. Then given any small $\sigma>0$,
there is a map $g\in{\tilde{\cal T}}$ such that
$$\mu<\delta (1+\sigma).  \eqno (5-9) $$
Let us start with such a map, and go through the general construction
described above. If this construction did not stop then we would have
a map of class ${\cal T}$ with arbitrary small initial scaling factor
such that the scaling factors of the renormalized maps would not go to 0.
Since this is impossible, our construction must stop. Then we come up
with a map $H=H_N$ of class $\tilde{\cal T}$ (with $T'=T^{N,j(N)}$ as the
central interval) such that $H T'\supset T'$ (high return). 
Since  the scaling factors of $H$ stay away from 0
(compare the above Remark on the scaling factors), we conclude that

$$\mu<{|T'|\over |T^0|} (1+\sigma)<{{\rm dist}(H c,c)\over |T^0|} (1+\sigma).
   \eqno (5-10)$$

\medskip\noindent
{\bf Push-forward.} Let us consider now the interval $B=T^{N,0}$ 
just constructed.
The map $H=f^{\circ(l+1)}$ is unimodal on $B$.
Moreover, there is an interval $A\supset f B$ which is a monotone  pull-back 
of $T^0$ by $f^l$. 

Let $p=(m,j)$ denote a point of the index set 
$$\tilde E=\{(m,j)\in E : m\leq N\quad {\rm or} \quad m=N, j=0\}.$$
We start with $0\equiv (0,0)$, and
by $p+1$ we mean the point of $\tilde E$ which lexicographically follows $p$.
(So, we identify the set $\tilde E$ ordered lexicographically
with an interval of the set of integers $p=0,...,P$).

Let us consider now a sequence of intervals 
$A^p=f^{\circ k(p)} A \subset T^{p-1}\quad p=,0,...P $
 defined as the last interval of
the orbit $\{ f^{\circ m} A\}_{m=0}^l$ visiting $T^{p-1}$,
$s(p)=k(p)-k(p+1)$.  Note that $A^0=T^0$ and $k(0)=l$, while
for $ p > 0,\; A^p \subset T^{p-1}\ssm T^p$. Moreover,
$f^{\circ s(p)}$ diffeomorphically maps $A^{p+1}$ onto $A^p$, and
the map $f^{\circ (s(p)-1)}| f A^{p+1}$ has a Koebe space
spread over $T^{p-1}$. Also the Koebe space of $f^{\circ k(P)} A$
is spread over $T^{P-1}$.  

Finally let us mark in $A^p$ the corresponding iterate
$a_p=f^{\circ(s(p)+1)}$ of the critical point. Then $f^{\circ s(p)}$
gives a diffeomorphism between corresponding marked intervals.
Note also that $a_0=H c$.

\medskip\noindent
{\bf Distortion estimates.}  For $1\leq p\leq P$ let the marked
point $a_p$ divide $A^p$ into intervals $L^p$ and $G^p$ with $G^p$ being
closer to the critical point than $L^p$. Set $\kappa_p=|G_p|/|L_p|$.

Somewhat abusing notations let us denote by $\mu_n=|T^n|/|T^{n-1}|$
new scaling factors.
Acting now as in the above particular cases taking into account estimate
(5-10) we obtain the following  analogue of (5-5):
$$\kappa_0\leq {1\over \mu_1} (1+o(1))\quad {\rm as}\quad \sigma\to 0. 
   \eqno (5-11)$$

Let us now estimate $\kappa_{p+1}$ through $\kappa_p$. We have (compare (5-6)):
$$ |G^p|\leq {\kappa_p\over 1+\kappa_p} |A^p|
   \leq {\kappa_p\over 1+\kappa_p}(1-\mu_p). \eqno (5-12) $$ 
 Let $R^p$ be the component of $T^{p-1}\ssm A^p$ containing $c$. Then
we can estimate the Poincar\'e length $P(L^p,R^p)$ as in (5-7):
$$P(L^p,R^p)\leq \log\left(1+{2\kappa_p\over 1+\mu_p}\right)\leq
    \log(1+{\kappa_p\over\mu_p}). \eqno (5-13)$$
By the Schwarz lemma
$${|f G^{p+1}|\over |f L^{p+1}|}\leq {\kappa_p\over \mu_p}. \eqno (5-14)$$
Since non-linearity of $f|A^{p+1}$ is estimated by $\log(1/\mu_{p+1})$,
$$ \kappa_{p+1}\leq {\kappa_p\over \mu_p\mu_{p+1}}. \eqno (5-15)$$ 
Estimates (5-11) and (5-15) yield
$$\kappa_p\leq {1\over (\mu_1...\mu_{P-1})^2\mu_P}(1+o(1)). \eqno  (5-16)$$
Set now $G'=A\ssm f B$. Then using (5-14)-like estimate and (5-16) we
conclude that
$${|G'|\over |f B|}\leq {1\over (\mu_1...\mu_{P})^2}(1+o(1)). \eqno  (5-17) $$
Finally, we can actually improve this estimate as
$${|G'|\over |f B|}\leq \tau^2 {1\over (\mu_1...\mu_{P})^2} =
    \left(\tau {|T^0|\over |B|}\right)^2\eqno  (5-18) $$
with an absolute $\tau<1$ (see the argument preceding estimate (5-8)).

\medskip\noindent
{\bf Complex pull-backs.} Take now the disk $D(T^0)$, and pull it back
by the branches of $H$. Then by (5-18) the central domain $V_0$ based
upon the interval $B$ is compactly contained in $D(T^0)$. The non-central
domains are compactly contained in $D(T^0)$, pairwise disjoint, and
disjoint from $V_0$ for the same reason as in the particular cases
treated above. 


 \bigskip
\centerline{\bf \S 6. Polynomial-like maps.}\medskip

The following
three pages, up to Lemma 6.3, present a self-contained exposition of
 a generalized version of the Branner-Hubbard
theory [BH], [B].  The reader can see the following differences.
We adjust the theory to a local setting of generalized  polynomial-like maps
which gives us a great flexibility in applications
 (see, e.g., [L2]) ( The original theory
was all about cubic polynomials with one escaping critical point).
We translate it from the original tableau language to the language of pull-backs
which nicely corresponds to the one-dimensional discussion of \S 2.  
 Finally,  we state the main rigidity result of
the theory in the parameter plane as a lemma on qc 
conjugacy of polynomial-like
maps with the same combinatorics (Lemma 6.3).
A direct proof of this lemma was given by J.Kahn. 

After this preparation we complete the proof of Theorem B.

\medskip\noindent{\bf Puzzle-pieces, kids and pull-backs.} 
Let us remember that
by a {\sl polynomial like map} we mean an analytic branched covering 
$$f: \cup_{i=0}^l U_i\rightarrow V$$ where $U_i$ and $V$ are topological disks,
cl$U_i\subset V$.  (The  Douady-Hubbard polynomial-like maps [DH]
correspond to $l=1$.) 
The set
$$K(f)=\{x: f^{\circ n} x\in\cup U_i,\; n=0,1,...\}$$
of non-escaping points is called the {\sl filled-in Julia set}.
By $\cal A$ we denote
a class  of \break polynomial-like maps
$f$ with a single non-degenerate critical point $c\in U_0$. This class
is a complex counterpart of the class ${\cal T}$ of one-dimensional maps.

\bigskip
\centerline{\psfig{figure=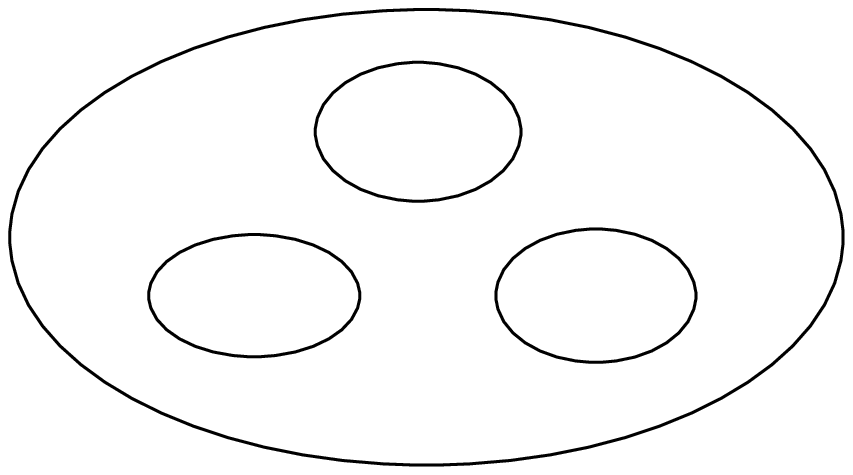,height=1.8in,width=4in}}\medskip
\centerline{Figure 7}\medskip
 
Set $V^{(n)}=f^{-n}V$.  The connected components $V^{(n)}_k$ of $V^{(n)}$
 are called
{\sl puzzle-pieces of depth } $n$.
 The puzzle piece of level $n$ containing a point $x$
will be also denoted by $V^{(n)}(x)$.  The puzzle-pieces 
$V^{(n)}_0\equiv V^{(n)}(c)$
 containing the critical point are
called {\sl critical}. The family of puzzle-pieces is Markov 
in the following sense:
for $n>0$ $V^{(n)}(x)$ is mapped under $f$ onto a puzzle piece 
$V^{(n-1)}(fx)$.
Moreover, this map is a two-to-one branched covering if 
$V^{(n)}(x)$ is critical,
 and a conformal isomorphism otherwise.

Let $f^{\circ m} x\in V^{(n)}_k\equiv W$. 
Then we can {\sl pull} the puzzle-piece $W$ {\sl back}
along the orb$_m(x)$ and  come up with the puzzle-piece 
$V^{(n+m)}(x)\equiv P$.
The pull-back is called {\sl univalent} if the map 
$f^{\circ m}: P\rightarrow W$ is. It is called {\sl quadratic-like}
if $P$ is critical and $f^{\circ(m-1)}: f P\rightarrow W$
is univalent.

A critical puzzle-piece $P=V^{(n+m)}_0$ is called a {\sl kid} of 
$W=V^{(n)}_0$ if
it is obtained by the quadratic-like pull back of $W$ along the
orb$_m(c)$.
If this corresponds to the  first return of
the critical point back to $W$, then $P$ is called
the {\sl first kid} of $W$. Repeating this construction we can talk about
{\sl grandkids} of the $n$th generation.

Let us say that $f$ {\sl admits a quadratic-like renormalization} if
there is a critical piece $W=V^n_0$ and its first kid $P=V^{n+m}_0$ such 
that the critical point does not escape $P$ under iterates of
$f^{\circ m}$. In such a case $f^{\circ m}: P\rightarrow W$ is a quadratic-like
map with a connected Julia set.

The critical point is called {\sl combinatorially recurrent} if its orbit
crosses all critical puzzle-pieces. 

\proclaim Two kids lemma. Assume that $f\in {\cal A}$ 
does not admit a quadratic-like
 renormalization. Then each critical puzzle-piece has at least two kids.

\noindent{\bf Proof.} Let us consider a critical puzzle-piece
 $W=V_0^{(n)}$ and its first kid $P=V_0^{(n+m)}$. Since $f$ does not admit
  quadratic-like renormalizations, the critical point must escape $P$
under some iterate $f^{\circ mk}$. Let $k$ be the first escape moment;
then $f^{\circ km} c\in W\ssm P$. Since the critical point is combinatorially
recurrent, we can find the first return moment $l$ of $f^{\circ km} c$ back
to $W$. Then the puzzle-piece $V^{(n-l)}(f^{\circ km} c)\subset W\ssm P$
is obtained by the univalent pull-back of $W$ along the 
${\rm orb}_l (f^{\circ km} c)$.
Pulling this piece further along the orb$_{km}(c)$ until it first hits
the critical point, we will find the second kid.  \QED

Let us consider now multiply-connected domains 
$A^{(n)}(x)=V^{(n)}(x)\ssm V^{(n+1)}$.
\break
The mod$(A^{(n)}(x)$) can be defined as the reciprocal of the extremal length
of the family of (non-connected) curves separating the outer boundary 
component of $A^{(n)}(x)$ from all inner ones. 

\proclaim The divergence property. Let $f\in {\cal A}$ 
does not admit a quadratic-like
   renormalization. Then for any $z\in K(f)$
$$\sum_{n=0}^\infty {\rm mod} (A^{(n)}(z)) =\infty.$$

\noindent{\bf Proof.} Argue as Branner \& Hubbard.
 Let us concentrate on the principle case when
the critical point is recurrent, and $z=c$. All critical pieces are 
descendents of $V=V^{(0)}_0$, and can be graded by generations.

 By the previous lemma
there are at least $2^n$ grandkids in $n$th generation. Since 
 mod($A^{(n)}_j)={\rm mod} (V)/2^n$ for any such grandkid, the total sum of
moduli over $n$th generation is at least mod$(V)$. Hence the total sum
of moduli over all descendents is $\infty$.   \QED 

\proclaim Corollary 6.1. A map $f\in {\cal A}$ does not admit a quadratic-like
 renormalization if and only if the filled-in Julia set is Cantor.

For this reason we also call maps which don't admit quadratic-like 
renormalizations  {\sl Cantor polynomial-like}. In this case $K(f)$
certainly coincides with the Julia set.

A set $K\subset{\bf C}$ is called {\sl removable} if given a neighborhood
$U\supset K$, any conformal/qc embedding $\psi: U\ssm K\rightarrow \bar{\bf C}$
allows conformal/qc continuation across $K$. The conformal and qc
settings in this definition are equivalent (by the Measurable Riemann Mapping
Theorem). The following important observation was made by Jeremy Kahn:

\proclaim Corollary 6.2.  Assume that $f\in {\cal A}$ does not admit a 
quadratic-like renormalization. Then the  Julia set $K(f)$ is removable. 

\noindent{\bf Proof.} This follows from the Modular Test on removability
                (see [SN], \S 1). \QED

 Talking about a conjugacy between two polynomial-like maps, we always mean
{\sl local} conjugacy in neighborhoods of their filled-in Julia sets.

\proclaim Lemma 6.3. Let $f$ and $g$ be two
 Cantor polynomial-like maps.
 If they are topologically conjugate by a homeomorphism $h$ 
then they are qc conjugate by a qc map $H$ which agrees with $h$ on the
Julia set.

\noindent {\bf Proof}. As in Lemma 2.10, let us set $g=\tilde f$ and mark the
related objects by tilde. Let $h$ be a topological conjugacy. Select an $N$
such that $V^{(N)}\subset {\rm Dom} (h)$. Let us consider an isotopy $h^t$
such that $h^0=h$, $h^1$ is smooth in a neighborhood of $V^{(N)}\ssm V^{(N+1)}$,
$ h^1\circ f=g\circ h^1 $ holds in a neighborhood of $\partial V^{(N+1)}$, 
and $h^t\equiv h$ in a neighborhood of the filled-in Julia set.
Since $h^t \equiv h$ near $K(f)$, we can pull this isotopy back to $V^{(N+1)}$
in such a way that for the pull-back $h_1^t\equiv h$ also holds near $K(f)$:

$$\matrix{      &     h^t_1     &             \cr
           \quad V^{(N+1)} & \rightarrow &\tilde V^{(N+1)}  \cr
          f\downarrow&       &\quad\downarrow \tilde f\cr
          \quad V^{(N)}        &  \rightarrow &\tilde V^{(N)}\cr
                   &     h^t        &        \cr}$$
   Then $h^t_1\equiv h^t$ in a neighborhood of $\partial V^{(N+1)}$. Hence we
can continue $h^t_1$ to $V^{(N)}$ as $h^t$.  

Now we can pull $h^t_1$ back in the same way, etc. We will obtain a
sequence of isotopies $h^t_n$ such that  

\noindent (i) $h_n^t$ agree with $h$ on $K(f)$;

\noindent (ii) $h_n^t$ agrees with $h^t_{n-1}$ on $V^{(N)}\ssm V^{(N+n)}$;

\noindent (iii) $\tilde f\circ h_n^t= h_{n-1}^t\circ f$.

\noindent (iv) $h_n^1$ is smooth outside $V^{(n)}$ with a uniformly bounded
  qc dilatation (since the pull-backs by conformal maps preserve
  the dilatation). 

Hence we can consider a family of maps $H^t=\lim h^t_n$ where the limit is
understood in a pointwise sense.  Then $H^t$ is an isotopy outside
the Julia set, and $H^1|V^{(N)}\ssm K(f)$ is a smooth qc map.

 Moreover, since
the isotopy $h^t_n$ is concentrated in $V^{(N+n)}$, it carries all puzzle-pieces
$V^{(N+n)}(x)$ to the corresponding pieces ${\tilde V}^{(N+n)}(hx)$. 
Hence  $H^1$ also carries $V^{(N+n)}(x)$  to  ${\tilde V}^{(N+n)}(hx)$. 
Since the diameters of these pieces shrink down to zero, we conclude
that $H^1$ is continuous across the Julia set, and agrees with $h$.

By Corollary 6.2, $H^1$ is actually  qc.  \QED

\smallskip
In conclusion let us mention the following result which links local and
global settings of the above theory (compare [DH]):

\proclaim Straightening Theorem. Any polynomial-like map of class ${\cal A}$
is  qc conjugate to a polynomial with one non-escaping critical point. \QED

\medskip\noindent {\bf R-symmetric case.}
If the topological disks $U_i$ and $V$ are symmetric with respect to
the real axis, and $f$ preserves the real axis, we call $f$ the 
${\bf R}$-{\sl symmetric polynomial-like map}. 

\proclaim Lemma 6.4. If two ${\bf R}$-symmetric polynomial-like 
maps $f$ and $g$ are
topologically conjugate on the real line by  
then the conjugacy can be continued to the complex plane as well.

\noindent {\bf Proof.} Let us continue the conjugating homeomorphism $h$
 to the complex plane in such
a way that it respects the dynamics on $\partial V^1$. Now let us pull it
back  (as in lemmas 2.10 and 6.30), 
so that the pull-backs $h_n: V\rightarrow V$ 
agree with $h$ on the real line. Then there is a pointwise limit $H=\lim h_n$
conjugating $f$ and $g$.
Since the puzzle-pieces shrink down to zero, $H$ is a homeomorphism.  \QED

\proclaim Lemma 6.5.
A {\bf R}-symmetric polynomial-like map $f\in {\cal A}$ 
 admits a quadratic-like renormalization if and
only if its restriction to the real line  admits a unimodal 
renormalization.

\noindent{\bf Proof.}
 Indeed, set $J=V\cap {\bf R}$ and 
$J^{(n)}_i=V^{(n)}_i\cap{\bf R}$ provided the
intersection is non-empty.  Then $J^{(n)}_i$ are exactly
the intervals of the Markov family ${\cal M}_J$ (see \S 2). 
Let $V^{(n+m)}_0$ be the first kid of $V^{(n)}_0$.
Then the  property ``orb($c$) does not escape $V^{(n+m)}_0$ under iterated
 $f^{\circ m}$"
is certainly equivalent to   
``orb($c$) does not escape $J^{(n+m)}_0$ under iterated $f^{\circ m}$".
The former property means that $f$ admits a quadratic-like renormalization,
while the latter one is equivalent to $f|{\bf R}$  admits a
unimodal renormalization (see Lemma 2.7).  \QED

Let us remember that ${\cal T}_{\min}$ denotes the class of maps 
$f\in {\cal T}$ with the recurrent $c$ and the
minimal critical set $\omega(c)$ which does not admit unimodal
renormalizations. Putting  together Theorem 2.12
and the last three lemmas, we conclude:

\proclaim Theorem 6.6. Two ${\bf R}$-symmetric polynomial-like maps 
 of class ${\cal T}_{\min}$ with the same combinatorial type
(that is $\bar\chi_f=\bar\chi_g$) are
qc conjugate.   \QED

Hence these maps are qs conjugate on the real line.\medskip

\noindent {\bf  The standard family.} Given  $p,q\in {\bf N}$ and a 
spin function $\epsilon$, let us  consider a standard
family of maps of class ${\cal T}(p,q,\epsilon)$
$$f^t:\bigcup_{k=-q}^p I_k\rightarrow J$$
defined after Lemma 2.7. The quadratic  central branch of $f^t$  depends
on $t$ while all non-central linear branches are fixed. By Lemma 2.7
any admissible combinatorial type $\bar\chi$ can be realised in such a
family. Since the lengths of the intervals $I^k$ can be selected arbitrarily,
we can construct a map $f$ with a given combinatorial type and
arbitrarily small first  scaling factor $\mu_0=|I_0|/|J|$.
 By \S 3, the conclusion of
Theorem B holds for such a map.

 Let us put the intervals $I_k$ into $J$ 
in such a way that they divide $J$ into commensurable parts. Then the pull-back
$U_0$ of the Eucledian disk $D(J)$ by the central branch has a bounded shape
(regardless  of the combinatorics and the lengths of $I_k$).
 All non-central pull-backs $U_k$ are true disks. Hence if $I_k$ are
sufficiently small, $f$ is polynomial-like.

So, we have constructed a polynomial-like map $f\in {\cal T}_{\min}$ with a
given combinatorial type satisfying the conclusion of Theorem B:
 Poincar\'{e} lengths of the gaps go up to $\infty$.
 But clearly this property
 is qs-invariant (since qs maps carry commensurable adjacent
intervals to commensurable ones). Now Theorem 6.5  yields:   

\proclaim Lemma 6.7.
 Theorem B holds for any polynomial-like map $f\in {\cal T}_{\min}$. \QED

\noindent{\bf Proof of Theorem B: concluding argument.} Let us consider the
subclass ${\cal T}^*$ of maps $f\in {\cal T}$ for which the conclusion
of Theorem B is valid.
 By the real argument of \S 3 and \S 4 this subclass includes all 
 maps with the non-minimal critical set,
 as well as maps with the minimal critical set of unbounded type.

Let us supply ${\cal T}$ with a $C^1$-topology. (Observe that
the  classes ${\cal T}(p,q,\epsilon) $ corresponding
to the different combinatorics on the first level (see \S 2)
stay far away even in
$C^0$-topology.) By  \S 3, ${\cal T}^*$ is an open subspace in ${\cal T}$.
Indeed, given an $\epsilon>0$, the condition that there is an $\epsilon$-small
 scaling factor $\mu_n$ specifies an open set of maps $f$. But by
\S 3 this condition  forces $f$ to belong to ${\cal T}^*$, provided 
$\epsilon $ is sufficiently small.

Let us take now a map $f\in {\cal T}_{\min}$ of bounded type. 
Assume that $f\not\in {\cal T}^*$. Then by \S 4
we can select a sequence of renormalized maps $R^{\circ n(k)} f$ 
$C^1$-converging to a map $g\in {\cal T}_{\min}$ of Epstein class.
By \S 5, there is a polynomial-like renormalization $h$ of $g$.
By Lemma 6.7, $h\in {\cal T}^*$, hence $g\in {\cal T}^*$.

Since ${\cal T}^*$ is open, $R^{\circ n(k)} f\in {\cal T}^*$ for sufficiently
large $k$. Hence $f\in {\cal T}^*$ as well, and this is a Contradiction.

 \noindent {\bf The case of a quadratic polynomial: alternative argument}.
 This case
can be treated in a more straightforward manner skipping 
\S\S 4,5. But then we need the Markov family of Yoccoz puzzle-pieces
(see [H] or [M2]) and the renormalization construction of [L2]. 

This construction goes as follows. Let $f$ be a non-tunable quadratic
polynomial in the sense of [DH]. (A real quadratic polynomial is
tunable if and only it admits a unimodal renormalization with period 
$>1$). Take a critical puzzle-piece $W=V^{(N)}_0$ 
such that $V^{(N-1)}_0\ssm V^{(N)}_0$ is a non-degenerate annulus.
 It satisfies the 
following ``nice" property similar to (2-4):
$$ f^{\circ n} (\partial W)\cap {\rm cl} W =\emptyset,\; n=1,2,...$$
If $\omega(c)$ is a minimal Cantor set then we can 
pre-renormalize $f$ on $W$
in the same way as it was described in \S 2 for the real setting.
Namely, consider the first return map to $W$ and select only those pieces
of its domain which intersect $\omega(c)$.   We obtain a polynomial-like
map $g:\cup U_i\rightarrow W$. It does not admit a quadratic-like
renormalization since $f$ is non-tunable. Hence it is a Cantor polynomial-like
map.

Now if we start with a real quadratic polynomial $f$ then Lemma 6.7 for $g$
implies Theorem B for $f$.

\bigskip
\centerline{\bf \S 7. Absence of attractors.}\medskip

Let $f$ be a quasi-quadratic map topologically exact on $T=[fc, f^{\circ 2}c]$.
By [BL1-3] or by [GJ], [K], $f|T$ has a unique measure-theoretic attractor
$A$, that is, an invariant closed set such that $\omega(x)=A$ for 
Lebesgue almost all $x\in T$. Moreover, either $A=T$ or $A=\omega(c)\ni c$.
Our goal is to prove that the former case holds. It is certainly true if
$c$ is not recurrent, or if $\omega(c)=T$. 

Let us assume that $c$ is recurrent and $\omega(c)\not= T$.
Then we can renormalize $f$ on any nice interval $I^0$.
 Moreover, for $I^0$ sufficiently short
 the domain of the pre-renormalized map
$g=g_1 :\cup I^1_k\rightarrow I^0$ is not dense, that is
$${\rm int} (I^0\ssm \cup I^1_k)\not= \emptyset. \eqno (7-0) $$
Since $g$ does not admit unimodal renormalizations,
the set $K(g)$ of non-escaping points is nowhere dense (see the argument
of Lemma 2.3).
Our goal is to prove that this set has zero Lebesgue measure.

\medskip\noindent{\bf Combinatorics of the first return maps.}
Let
$I^0\supset I^1\supset I^2\supset...$ be the sequence of the first kids
of $I^0$.
 Let us construct inductively the return maps $f_n$ to $I^{n-1}$. Let 
${\bf L}^n=\cup L^n_i$ be the domain of definition of the $f_n$
where $L^n_i$ are intervals with $L^n_0$ as the central interval,
${\bf L}^n_*=\cup_{i\neq 0} L^n_i$.
 Denote by ${\cal L}^n$ the family of these
intervals, and by ${\cal L}_*^n$ the family of non-central intervals.
Let $G^n={\rm int} (I^{n-1}\ssm {\bf L}^n)$ be the union of gaps.
 All these sets are $c$-symmetric for $n>1$.

To start induction observe that ${\bf L}^1 = \cup I^1_k,\; f_1 = g_1$.
Moreover, $G^1$ is non-empty by (7-0).
Let us denote by ${\cal G}^n_*$ the free semigroup generated by ${\cal L}_*^n$.
For a word $\gamma\in {\cal G}^n_*$ denote by the same letter 
$\gamma$ the interval
whose itinerary through ${\bf L}^n$ is given by $\gamma$. Let $l$ be the length
of $\gamma$. Let us
consider the subsets  $\gamma_r$ and $\gamma_e$ of this interval
such that $f_n^{\circ l}\; \gamma_0^n=I^n$ and 
 $f_n^{\circ l} \;\gamma_e=G^n$, that is, the former subset goes
to the central interval, while the latter one escapes through the
gaps set $G^n$. Finally, let $\partial {\cal G}^n$ mean the set of infinite
words in letters ${\cal L}^n_*$ which we identify with the corresponding
$f_n$-invariant Cantor set. Since this set does not contain the critical
point, it is hyperbolic and has zero measure. Hence 
the sets $\cup\gamma_e$ and $\cup\gamma_r$ cover almost completely the set
${\bf L}^n_*$. Let 
$$E^n=G^n\bigcup_{\gamma\in {\cal G}^n_*}\gamma_e,
\quad   R^n=I^n \bigcup_{\gamma\in {\cal G}^n_*}\gamma_r$$
be the full sets of returning and escaping points. Because
of the above remark, their union has full measure in $I^{n-1}$.
Moreover, we have a {\sl transition map} $F_n$ from $R^n$ to the central
interval $I^n$ which maps diffeomorphically each interval $\gamma_e$ onto
$I^n$.

Now, in order to construct the return map $f_{n+1}$ of the next level,
just consider the pull-back ${\bf L}^{n+1}$ of $E^n$ by the central
 branch $f_n|I^n$, and define $f_{n+1}=F_n\circ f_n|\; {\bf L}_{n+1}$.

\eject
\medskip\noindent{\bf Density estimates.}

\proclaim Lemma 7.1. Let $n>1$. Then for any word $\gamma\in {\cal G}^n_*$
$${|\gamma_e|\over |\gamma_r|}\geq  {1\over 2}{|G^n|\over |I^n|}.$$

\noindent {\bf Proof.} Let $l$ be the length of the word $\gamma$. Then
$f_n^{\circ l}$ diffeomorphically maps $\gamma$ onto $I^{n-1}$. By the Minimum
Principle, there is a component $L$ of $I^{n-1}\ssm I^n$ such that
$$\inf_L Dg_n^{\circ(-l)}\geq \sup_{I^n} Dg_n^{\circ(-l)}.$$ 
Let $H=G^n\cap L$. By the symmetry, $|H|=(1/2) |G^n|$. Hence
$${|\gamma_e|\over |\gamma_r|}\geq 
{|g_n^{\circ(-l)} H|\over |g_n^{\circ(-l)} I^n|}\geq {|H|\over |I^n|}=
{1\over 2}{|G^n|\over |I^n|}.$$  \QED

\proclaim Lemma 7.2. Let $U$ be any interval such that 
$\partial U\subset \partial {\bf L}^n\cup\partial I^{n-1}$
   and $U\cap I^n=\emptyset$. Then
$${{\rm dens} (E^n|U)\over {\rm dens} (R^n|U)}\geq 
{1\over 2}{|G^n|\over |I^n|}.$$  

\noindent {\bf Proof.} There is a family ${\cal L}$ of  intervals
$\gamma\in {\cal G}^n_*$ such that we have the following coverings up to
sets of measure zero:
$$(U\cap E^n)\ssm G^n=\bigcup_{\gamma\in\cal L}\gamma_e\; ({\rm mod}\; 0),\quad 
U\cap R^n=\bigcup_{\gamma\in\cal L}\gamma_r\; ({\rm mod}\; 0).$$
Apply now the previous lemma.     \QED

Let us state an elementary lemma about the quadratic map
$\phi: x\mapsto (x-c)^2$.

\proclaim Lemma 7.3. Let $X$ be a $c$-symmetric measurable set in a 
$c$-symmetric interval $I$. Then 
$${\rm dens} (X|I)\geq {1\over 2} {\rm dens} (\phi X|\phi I). $$\QED

 Let us remember the notation $\mu_n=|I^n|/|I^{n-1}|$ for the scaling
factors. Let $\delta_n={\rm dens}(G^n|I^{n-1}).$ If $U$ is  any interval
such that $\partial U\subset \partial {\bf L}^n\cup\partial I^{n-1}$
(but perhaps $U\supset I^n$) then Lemma 7.2 implies
$${{\rm dens} (E^n|U)\over {\rm dens} (R^n|U)}\geq 
  {1\over 2} \delta_n {1-\mu_n\over \mu_n}. \eqno (7-1) $$
Pulling this back by $f_n|I^n= h\circ\phi$ with the $h$ of bounded
distortion we obtain the following recurrent estimate:
$$\delta_{n+1}\geq a \delta_n {1-\mu_n\over \mu_n} \eqno (7-2) $$
with an absolute $a$. It follows that $\delta_{n+1}\geq 2\delta_n$,
provided $\mu_n\leq \bar\mu$ is sufficiently small.

If $\mu_m>\bar\mu$ is not too small then by the results of \S 3 
we are in the tail of a cascade of central returns. Let $m+q-1$
be the first non-central return level of this cascade, $q\geq 1$. 
Let $H$ be the component
of $I^{m+q-1}\ssm I^{m+q}$ containing the critical value $f_{m+q} c$.
Let us consider the intervals $V=H\cap f_{m+q} I^{m+q}$
and $U=V\cup I^{m+q}_k$ where $f_{m+q} c\in I^{m+q}_k$.
Then $U$ satisfies the assumptions of Lemma 7.2. On the other hand,
by \S 3, the Poincar\'{e} length $P(I^{m+q}_k | H)$ is very small.
We conclude that 
$${\rm dens} (G^{m+q}|V)\geq (1-\epsilon)\; 
{\rm dens} (G^{m+q}|U) \eqno (7-3)$$
with a very small $\epsilon$.

Further, $f^q: U\rightarrow f^q U$ is a map of bounded distortion.
Indeed, since $|I^{m+q}|/ |I^m|\geq \tau \bar\mu, \;\tau>0,$ stays
away from 0, $f_m$ has a bounded non-linearity on the interval
$I^m\ssm I^{m+q}$. Since the intervals $f^k U,\; k=0,...,m,$ are
pairwise disjoint, we have the bounded distortion of $f^q$.
Hence there is an absolute $a>0$ such that
$${\rm dens} (G^{m+q}|U)\geq a\; {\rm dens} (G^m| f^q U)  \eqno (7-4) $$
 Now Lemma 7.2 and estimates (7-3), (7-4) yield
$${\rm dens} (G^{m+q}|V)\geq a \delta_m \eqno (7-5)$$
(as usually in analysis, an absolute constant $a$ may have different
values in different estimates).

Now let us pass to the next level using Lemma 7.3. Let $W\subset I^{m+q}$ 
be the pull back of $V$ by $f_{m+q}|I^{m+q}$. We conclude that
$${\rm dens} (G^{m+q+1}|W)\geq a\delta_m \eqno (7-6)$$
with yet another $a>0$. On the other hand, by \S 3 
$|I^{m+q+1}|/| W|$ is very small (exponentially small in terms of
$\kappa$). Hence
$${\delta_{m+q+1}\over\mu_{m+q+1}}\geq {|G^{m+q+1}|\over |I^{m+q+1}|}\geq 
  {{\rm dens} (G^{m+q+1}|W)\over {\rm dens} (I^{m+q+1}|W)}>
   A\delta_m$$
with a big $A$ (for a big $m$).
 Passing now to the next level using (7-2), we conclude

\proclaim Theorem 7.4. The densities $\delta_n={\rm dens} (G^n|I^{n-1})$
of gap sets stay away from 0, provided
$n$ is not in the tail of a long cascade of central returns or 
immediately after the cascade. Moreover, these densities
  grow at least exponentially with $\kappa(n)$.

\noindent {\bf Concluding argument.} Assume that the set $K(g)$
of non-escaping points has positive measure. Let 
$X=\{x\in K(g): \omega(x)\ni c \}$. By [BL3],
dens$(X|I^n)=1$ (The argument:
 take a density point $x\in X$ and consider the first moment $l$
when $g^l x\in I^{n}$. Then the corresponding pull-back
 $T\ni x$ is
mapped under $g^l$ onto $I^n$   with a bounded distortion.)

On the other hand, Theorem 7.4 says that dens$(K(g)|I^n)\to 0$
for an appropriate subsequence of levels.
This contradiction completes the proof of theorem A.

\eject

\centerline{\bf References.}\smallskip

{\item [B]}. B. Branner. Cubic polynomials: turning around the
  connectedness locus. Preprint of  The Technical University of
  Denmark, 1992-05. 
  To appear in
    ``Topological Methods in Modern Mathematics, A Symposium in
     Honor of John Milnor's 60th Birthday". \smallskip

{\item [BH]} B.Branner \& J.H.Hubbard. The iteration of cubic polynomials,
   Part II   : patterns and parapatterns, Acta Math., to appear.\smallskip

{\item[BL1]} A.M.Blokh \& M.Lyubich.
 Attractors of transformations  of an interval.
 Functional Analysis and Applications, {\bf 21} (1987), $n^\circ 2$, 148-150.
  \smallskip

{\item[BL2]}  A.M.Blokh \& M.Lyubich. Typical  behavior of the trajectories 
   of transformations of a segment. Teoriya Funktsii, Funktsional'nyi
  Analiz i ikh Pril., {\bf 49} (1988), 5-16. Translated in Journal of
   Soviet Math, {\bf 49} (1990), 1037-1044. \smallskip

{\item [BL3]} A.Blokh \& M.Lyubich. Measurable dynamics of S-unimodal maps
    of the interval. Ann. scient. \'{E}c. Norm. Sup., v. 24 (1991),
    545-573. 

{\item [DH]} A.Douady \& J.H.Hubbard. On the dynamics of polynomial-like maps,
    Ann. scient. \'{E}c. Norm. Sup. v. 18 (1985), 287-343.\smallskip

{\item [G]} J.Guckenheimer. Sensitive dependence to initial conditions 
     for one-dimensional maps. Comm. Math. Phys., v.70 (1979), 133-160.
     \smallskip 

{\item [GJ]} J. Guckenheimer \& S. Johnson. Distortion of $S$-unimodal maps.
        Annals Math., v. 132 (1990), 71-130.\smallskip

{\item [JS]} M.Jacobson \& G.\'{S}wiatek. Metric properties of 
     non-renormalizable $S$-unimodal maps. I. Preprint IHES/M/91/16.

{\item [JR]} L. Jonker, D. Rand. Bifurcations in one dimension. I. The 
    non-wandering set. Inventions Math., v.62 (1981), 347-365.\smallskip

{\item [H]} J.H.Hubbard. Local connectivity of Julia sets and bifurcation
     loci: three theorems of J.-C. Yoccoz. To appear in
    ``Topological Methods in Modern Mathematics, A Symposium in
     Honor of John Milnor's 60th Birthday".\smallskip

{\item [HK]} F. Hofbauer and G. Keller. Some remarks on recent results 
         about S-unimodal maps. Ann. Inst. Henri Poincar\'{e}, v. 53, 
         \#4 (1990), 413-425.\smallskip

{\item[K]} G.Keller. Exponents, attractors and Hopf decomposition
  for interval maps. Ergodic Theory \& Dynamical Systems, {\bf 10} (1990),
   717-744.\smallskip

{\item [KN]} G. Keller \& T. Nowicki. Fibonacci maps revisited.
             Preprint, 1992. \smallskip
 
{\item [L1]} M. Lyubich.  Non-existence of wandering intervals and
    structure of topological attractors for one-dimensional dynamical systems.
    Erg. Th. \& Dyn Syst. {\bf 9} (1989), 737-750.\smallskip

{\item [L2]} M.Lyubich. On the Lebesgue measure of the Julia set of a 
     quadratic polynomial. Preprint IMS Stony Brook, 1991/10. \smallskip

{\item [L3]} M. Lyubich. Milnor's attractors, persistent recurrence and
  renormalization. To appear in
    ``Topological Methods in Modern Mathematics, A Symposium in
     Honor of John Milnor's 60th Birthday".\smallskip

{\item [LM]} M.Lyubich \& J.Milnor. The unimodal Fibonacci map. 
        Preprint \#15/1991, Stony Brook.\smallskip 

{\item [M1]} J.Milnor. On the concept of attractor. Comm. Math. Phys,
     {\bf 99} (1985), 177-195, and {\bf 102} (1985), 517-519. \smallskip

{\item [M2]} J.Milnor. Local connectivity of Julia sets: expository
   lectures. Preprint IMS Stony Brook, \#1992/11.\smallskip
 
{\item [MT]} J.Milnor \& W.Thurston. On iterated maps of the interval,
           pp. 465-563 of ``Dynamical Systems, Proc. U. Md., 1986-87,
          ed. J. Alexander, Lect. Notes Math., 1342, Springer 1988.\smallskip

{\item [Ma]} M.Martens. Distortion results and invariant Cantor sets
          of unimodal maps. Preprint IMS Stony Brook, \# 1992/1.
          (A part of the thesis ''Interval dynamics", 1990).
              \smallskip

{\item [MS]} W. de Melo \& S. van Strien. One dimensional dynamics. \smallskip

{\item [NS]} L. Sario \& M. Nakai. Classification Theory of Riemann Surfaces.
              Springer-Verlag, 1970.\smallskip
 
{\item [S]} D.Sullivan. Bounds, quadratic differentials, and renormalization
         conjectures,  1990. To appear in AMS Centennial Publications.
        {\bf 2}: Mathematics into Twenty-first Century. \smallskip

{\item [vS]} S. van Strien. On the bifurcations creating horseshoes.
           Springer Lect. Notes Math., v. 898  (1981), 316-351.\smallskip

\bye